\documentclass[review,3p,number]{elsarticle}
\usepackage[dvipsnames, table]{xcolor}
\usepackage{amsmath}
\usepackage{amssymb}
\usepackage{amsthm}
\usepackage[font = footnotesize]{caption}
\usepackage{float}
\usepackage{hyperref}
\usepackage{graphicx}
\usepackage{natbib}
\usepackage{nicematrix}
\usepackage{xparse}

\DeclareMathOperator* {\supp} {supp}

\newtheorem{thm}{Theorem}[section]
\newtheorem{lemma}[thm]{Lemma}

\newtheorem{prop}[thm]{Proposition}
\newdefinition{remark}{Remark}[section]
\NewDocumentEnvironment{pf} {o}
{
    \noindent
    \IfValueTF{#1} {\textit{Proof of #1.}} {\textit{Proof.}}
}
{\hfill $\square$ \vspace{6pt}}

\newcommand {\dif} [1] {\, \text{d} #1}

\DeclareDocumentCommand {\dd} {m o} 
    { \IfValueTF{#2}
        { \left( \frac{\text{d}} {\text{d} #1} \right)^{#2} \! } 
        { \frac{\text{d}} {\text{d} #1} } 
    }
\DeclareDocumentCommand {\complex} {o} { \mathbb{C}\IfValueT{#1}{^{#1}} }
\DeclareDocumentCommand {\real}    {o} { \mathbb{R}\IfValueT{#1}{^{#1}} }

\hypersetup {
    colorlinks = true,
    linkcolor = black,
    filecolor = magenta,
    urlcolor = NavyBlue,
    citecolor = NavyBlue
}
\definecolor{myblue} {HTML}{34CDF9}
\definecolor{mygreen}{HTML}{67FD9A}
\definecolor{mygrey} {HTML}{9B9B9B}
\definecolor{myred}  {HTML}{FF8F8C}

\begin{document}
\begin{frontmatter}
    \title{Uniqueness of the Inverse Source Problem for Fractional Diffusion-Wave Equations}

    \author[Qiu1,Qiu2]{Lingyun Qiu} 
    \ead{lyqiu@tsinghua.edu.cn}

    \affiliation[Qiu1]{
        organization={Yau Mathematical Sciences Center, Tsinghua University},
        city={Beijing},
        postcode={100084},
        country={China}
    }
    
    \affiliation[Qiu2]{
        organization={Yanqi Lake Beijing Institute of Mathematical Sciences and Applications},
        addressline={No. 544, Hefangkou Village, Huaibei Town}, 
        city={Huairou District, Beijing},
        postcode={101408}, 
        country={China}
    }
    
    \author[Sim]{Jiwoon Sim\corref{cor}\fnref{adr}}
    \cortext[cor]{Corresponding author.}
    \ead{jiwoonsim@gmail.com}
    \fntext[adr]{Present address: Department of Mathematical and Statistical Sciences, University of Alberta, Edmonton, T6G2R3, Alberta, Canada.}

    \affiliation[Sim]{
        organization={Department of Mathematical Sciences, Tsinghua University},
        city={Beijing},
        postcode={100084}, 
        country={China}
    }

    \begin{abstract}
        This study addresses the inverse source problem for the fractional diffusion-wave equation, characterized by a source comprising spatial and temporal components. The investigation is primarily concerned with practical scenarios where data is collected subsequent to an incident. We establish the uniqueness of either the spatial or the temporal component of the source, provided that the temporal component exhibits an asymptotic expansion at infinity. Taking anomalous diffusion as a typical example, we gather the asymptotic behavior of one of the following quantities: the concentration on partial interior region or at a point inside the region, or the flux on partial boundary or at a point on the boundary. The proof is based on the asymptotic expansion of the solution to the fractional diffusion-wave equation. Notably, our approach does not rely on the conventional vanishing conditions for the source components. We also observe that the extent of uniqueness is dependent on the fractional order.
    \end{abstract}



    \begin{keyword}
        Fractional diffusion-wave equation (primary) \sep Inverse source problem \sep Asymptotic expansion \sep Mittag-Leffler functions

    \end{keyword}
\end{frontmatter}
    
    \section{Introduction}
    \label{sec: intro}
    Let $\Omega \in \real[d]$ be a bounded region with a smooth boundary $\partial \Omega$. Let $A$ be a uniformly elliptic second-order differential operator with a nonnegative potential, expressed as
    \begin{equation*}
        Au(x) = -\sum_{i, j = 1}^d \partial_i (a_{ij}(x) \partial_j u(x)) + q(x) u(x)
    \end{equation*}
    where $a_{ij} = a_{ji} \in C^1(\overline{\Omega})$, $q \in C(\overline{\Omega})$ and $q \geq 0$. Moreover, there exists $c_0 > 0$ such that
    \begin{equation*}
        \sum_{i, j = 1}^d a_{ij}(x) \xi_i \xi_j \geq c_0 ||\xi||^2
    \end{equation*}
    for all $x \in \overline{\Omega}$ and $\xi = (\xi_i)_{1 \leq i \leq d} \in \real[d]$. Let us consider the following initial boundary value problem for a fractional diffusion-wave equation:
    \begin{equation}
        \left\{
        \begin{aligned}
            \partial_t^\alpha u + Au &= f(x) \mu(t) && \text{in } \Omega \times (0, \infty) \\
            u(0)                     &= 0           && \text{in } \Omega \ (0 < \alpha \leq 2) \\
            \partial_t u(0)          &= 0           && \text{in } \Omega \ (1 < \alpha \leq 2) \\
            u                        &= 0           && \text{on } (\partial \Omega) \times (0, \infty)
        \end{aligned}
        \right. \label{eq: TFDE} 
    \end{equation}
    Here $\alpha \in (0, 1) \cup (1, 2)$, $\partial_t^\alpha$ is the Caputo fractional derivative
    \begin{equation*}
        \partial_t^\alpha g(t) = \left\{
        \begin{aligned}
            & \frac{1}{\Gamma(k - \alpha)} \int_0^t g^{(k)}(\tau) (t - \tau)^{k - 1 - \alpha} \dif{\tau}, && k - 1 < \alpha < k \\
            & \partial^k g(t),  && \alpha = k \in \mathbb{N}
        \end{aligned}
        \right.
    \end{equation*}
    and $f$ and $\mu$ represent the spatial and temporal distribution of the source respectively. 
    
    When $\alpha = 1$ or $2$ in equation \eqref{eq: TFDE}, it becomes a regular diffusion or wave equation. Hence, a fractional diffusion-wave equation is a generalization of diffusive processes. Metzler and Klafter \cite{mk00} deduced this equation for $\alpha < 1$ by studying the stochastic movement of particles. The equation is shown effective to model various phenomena including diffusion in fractal or porous media \cite{n86, gvh06} and dielectric relaxation \cite{n84}. 
    
    The forward problem, solving $u$ from initial and boundary conditions, has been extensively studied in the literature. Sakamoto and Yamamoto \cite{sy11a} proposed a weak solution to equation \eqref{eq: TFDE}, while Gorenflo et al. \cite{gly15} relaxed the regularity of the source term. The maximum principle for classical and weak solutions was derived by Luchko \cite{l09} and Liu et al. \cite{lry16}.
    
    The inverse source problem seeks to determine $f$ and/or $\mu$ from some measurement of $u$. Inverse source problems can roughly be divided into inverse $x$-source problems and inverse $t$-source problems depending on which component to recover. We define the normal derivative with respect to $A$ by
    \begin{equation*}
        \partial_\nu u = \mathbb{A} \nabla u \cdot \nu,
    \end{equation*}
    where $\mathbb{A} = (a_{ij})_{1 \leq i, j \leq n}$ and $\nu$ is the normal derivative on $\partial \Omega$. Let $I$ be a non-empty time interval, $\omega$ be a non-empty subdomain of $\Omega$ and $\gamma$ be a non-empty open subboundary of $\partial \Omega$. In an inverse $x$-source problem, one recovers $f$ provided that $\mu$ is known and one of the following data is given: (a) $u|_{\omega \times I}$, (b) $\partial_\nu u|_{\gamma \times I}$ or (c) $u(\cdot, T)$ for some $T > 0$. Similarly, in an inverse $t$-source problem, one recovers $\mu$ provided that $f$ is known and one of the following data is given: (a) $u|_{\omega \times I}$, (b) $\partial_\nu u|_{\gamma \times I}$ or (c$_1$) $u(x_0, \cdot)$, (c$_2$) $\partial_\nu u(x_0, \cdot)$ depending on whether $x_0 \in \Omega$ or $x_0 \in \partial \Omega$. 
    
    In practical situations, one can only start observation after an incident has already happened. Hence, it is desirable that the support of $\mu$ includes 0 and the time interval $I$ has a positive distance from $0$. In recent research the uniqueness of the inverse problem under this constraint is proved with some additional assumptions. Yamamoto \cite{y23}, and Janno and Kian \cite{jk22}, independently proved the uniqueness of both inverse $x$- and $t$-source problems with data types (a) or (b). The observation region $\omega$ or $\gamma$, as well as the time interval $I$, can be arbitrarily small, but they assumed that the source stops for some time before the end of observation. Janno and Kian additionally addressed the inverse $t$-source problem with data type (c) and the simultaneous recovery of $f$ and $\mu$. A similar result not requiring the vanishing of the source before the end of observation was proved by Kian et al. \cite{kly22}. Nevertheless, a stronger assumption is imposed, in which the spatial component of the source should be zero in the observation area $\omega$. 

    Our work extends the aformentioned results to include the uniqueness of inverse problems where data is asymptotic and source vanishing conditions are not assumed. \\
    \noindent (IP1) Determine $f$ and/or $\mu$ from the asymptotic behavior of $u(\cdot, t) |_\omega$ or $\partial_\nu u(\cdot, t) |_\gamma$ as $t \rightarrow \infty$.\\
    \noindent (IP2) Determine $\mu$ from the asymptotic behavior of $u(x_0, t)$, $x_0 \in \Omega$ or $\partial_\nu u(x_0, t)$, $x_0 \in \partial \Omega$ as $t \rightarrow \infty$. \\
    In these inverse problems we neither assume the spatial component $f$ nor the temporal component $\mu$ vanishes. Moreover, since we only take asymptotic behaviors, the observation may not start at the beginning of an incident. This is of great significance to practical application, for example, the determination of the pollutant source beneath the soil. Since we do not a priori know the location of the source, we can not assume that $f$ vanishes in the observation area. On the other hand, the pollutant may be emitted continuously, and we do not know when or whether it will stop, so we can not assume the vanishing of $\mu$ either. To the best of our knowledge, this is the first result on the uniqueness of inverse source problems for fractional diffusion-wave equations such that the data is taken away from the beginning time and no vanishing conditions are imposed on the source.
    
    Let us briefly mention some other related works on inverse source problems. The uniqueness of inverse $x$-source problem with data type (c) is proved by Sakamoto and Yamamoto \cite{sy11b} and Slodi\v{c}ka \cite{s20}. Some stability results are given in \cite{lz17, kst23, kj18, cl23}. The inverse source problems for heat and wave equations are studied in the works of Jin Cheng et al. \cite{cly23} and Kian et al. \cite{kly22}. Readers can also refer to the review \cite{lly19} on inverse source problems and the references therein. For the convenience of readers we present an illustration of existing literature, our contributions and open problems in Table~\ref{table: summary}.

    \begin{table}[h!]
        \centering
        \begin{NiceTabular}{cc*{5}{Wc{36pt}}}[hvlines]
            \hline
            \Block{1-2}{ \diagbox{\vspace{-5pt} \hspace{12pt} $t$} {\vspace{2pt} $x \qquad$} } 
            & & $\Omega$ & $\omega \subseteq \Omega$ & $\gamma \subseteq \partial \Omega$ & $x_0 \in \Omega$ & $x_0 \in \partial \Omega$ \\ \hline
            \Block{3-1}{$T_1 = 0$} & $T_* < T_2 < \infty$ & \cellcolor{mygrey} & \cellcolor{mygrey}  & \cellcolor{mygrey} & \cellcolor{myred} & \cellcolor{myred} \\ \cline{2-7} 
            & $T_2 = \infty$ & \cellcolor{mygrey} & \cellcolor{mygrey} & \cellcolor{mygrey} & \cellcolor{myred} & \cellcolor{myred} \\ \cline{2-7} 
            & $T_2 < \infty$ & \cellcolor{mygrey} & \cellcolor{mygreen}\textit{\cite{kly22}} & \cellcolor{mygreen}\cite{cl23} & \cellcolor{myred} & \cellcolor{myred} \\ \hline
            \Block{3-1}{$T_1 > 0$} & $T_* < T_2 < \infty$ & \cellcolor{mygrey} & \cellcolor{mygreen}\cite{y23}   & \cellcolor{mygreen}\cite{y23} & \cellcolor{myred} & \cellcolor{myred} \\ \cline{2-7} 
            & $T_2 = \infty$ & \cellcolor{mygrey} & \cellcolor{myblue} & \cellcolor{myblue} & \cellcolor{myred} & \cellcolor{myred} \\ \cline{2-7} 
            & $T_2 < \infty$ & & \cellcolor{mygreen}\cite{kly22} & & \cellcolor{myred} & \cellcolor{myred} \\ \hline
            \Block{1-2}{$T_1 = T_2 = T$}  & & \cellcolor{mygreen}\cite{sy11b} \cite{s20} & & & \cellcolor{myred} & \cellcolor{myred} \\ \hline
        \end{NiceTabular}

        \vspace{12pt}
        \begin{NiceTabular}{cc*{5}{Wc{36pt}}}[hvlines]
            \hline
            \Block{1-2}{ \diagbox{\vspace{-5pt} \hspace{12pt} $t$} {\vspace{2pt} $x \qquad$} }  
            & & $\Omega$ & $\omega \subseteq \Omega$ & $\gamma \subseteq \partial \Omega$ & $x_0 \in \Omega$ & $x_0 \in \partial \Omega$ \\ \hline
            \Block{3-1}{$T_1 = 0$} & $T_* < T_2 < \infty$ & \cellcolor{mygrey} & \cellcolor{mygrey} & \cellcolor{mygrey} & \cellcolor{mygrey} & \cellcolor{mygrey} \\ \cline{2-7} 
            & $T_2 = \infty$ & \cellcolor{mygrey} & \cellcolor{mygrey} & \cellcolor{mygrey} & \cellcolor{mygrey} & \cellcolor{mygrey} \\ \cline{2-7} 
            & $T_2 < \infty$ & \cellcolor{mygrey} & \cellcolor{mygreen}\textit{\cite{kly22}} & & \cellcolor{mygreen}\cite{lz17}  & \\ \hline
            \Block{3-1}{$T_1 > 0$} & $T_* < T_2 < \infty$ & \cellcolor{mygrey} & \cellcolor{mygreen}\cite{y23} & \cellcolor{mygreen}\cite{y23} \cite{jk22} & \cellcolor{mygreen}\textit{\cite{y23}} & \cellcolor{mygreen}\textit{\cite{y23}} \cite{jk22} \\ \cline{2-7} 
            & $T_2 = \infty$ & \cellcolor{mygrey} & \cellcolor{myblue} & \cellcolor{myblue} & \cellcolor{myblue} & \cellcolor{myblue} \\ \cline{2-7} 
            & $T_2 < \infty$ & & \cellcolor{mygreen}\cite{kly22} & & & \\ \hline
            \Block{1-2}{$T_1 = T_2 = T$} & & \cellcolor{myred}  & \cellcolor{myred} & \cellcolor{myred} & \cellcolor{myred} & \cellcolor{myred} \\ \hline
        \end{NiceTabular}
        
        \caption{A summary of inverse $x$- (above) and $t$- (below) problems. The first row and first two columns show the region and time interval where data is gathered. $(T_1, T_2)$ is the time interval and $T_*$ is the stopping time of the source (if applicable). Cells with an \textit{italic} citation means that the case can be solved using the method in the article cited, with slight modifications in the details. Colour of the cell: grey --- can be covered by other existing results; red --- incompatible data; green --- resolved; blue --- our contribution.}
        \label{table: summary}
    \end{table}
    
    The rest of the article is organized as follows. In Section~\ref{sec: preliminaries}, we provide the preliminaries and present the main theorems. Section~\ref{sec: proof} is dedicated to the proofs of these theorems. Finally, Section~\ref{sec: conclusions} offers concluding remarks and discusses potential areas for future research.
    
    \section{Preliminaries and main theorems}
    \label{sec: preliminaries}
    This section extends the discussion to equation \eqref{eq: TFDE}. Define the domain of $A$ by $\mathcal{D}(A) = H^2(\Omega) \cap H_0^1(\Omega)$. Suppose $f \in L^2(\Omega)$ and $\mu \in L^\infty (0, \infty)$, then equation \eqref{eq: TFDE} admits a weak solution $u \in L^2_\text{loc} (0, \infty;$ $\mathcal{D}(A)) \cap H^\alpha_\text{loc} (0, \infty; L^2(\Omega))$ as established in \cite[Theorem 2.2]{sy11a}. Given that the spectrum of $A$ consists of positive eigenvalues of finite multiplicities, these can be denoted as $\{ \lambda_n \}_{n \in \mathbb{N}}$ in increasing order, without counting multiplicities.
    Let $P_n$ be the orthogonal projection onto the eigenspace of $\lambda_n$. The solution $u$ is given by (refer to \cite[Lemma 1]{y23})
    \begin{equation}
        u(x, t) = \sum_{n = 1}^\infty \left( \int_0^t (t - s)^{\alpha - 1} E_{\alpha, \alpha} (-\lambda_n (t - s)^\alpha) \mu(s) \dif{s} \right) P_n f(x),
        \label{eq: solution u}
    \end{equation}
    where $E_{\alpha, \beta}(z) = \sum_{k = 0}^\infty z^n / \Gamma(\alpha k + \beta), \, \alpha, \beta > 0$ is the Mittag-Leffler function. We will frequently use the following properties of the Mittag-Leffler function.
    
    \begin{prop}
        \label{prop: Mittag-Leffler}
        The Mittag-Leffler function $E_{\alpha, \beta}$ satisfies the following properties: \\
        (a) $E_{\alpha, \beta}$ is entire in $\mathbb{C}$. \\
        (b) $E_{\alpha, \beta}$ admits an asymptotic expansion (\cite[Theorem~1.4]{p99}): 
        \begin{equation*}
            E_{\alpha, \beta}(-x) \sim \sum_{k = 1}^\infty \frac{ (-1)^{k + 1} }{ \Gamma(\beta - k \alpha) } x^{-k} \text{ as } x \rightarrow +\infty.
        \end{equation*}
        Here $\frac{1}{\Gamma(z)} = 0$ for $z = 0, -1, -2, \ldots$. \\
        (c) $E_{\alpha, \beta}(-x) \leq C(1 + x)^{-1}$ for any $x \geq 0$ (\cite[Theorem~1.6]{p99}).
    \end{prop}
    
    For all $\sigma > 0$, we define the operator $A^\sigma$ by
    \begin{equation*}
        A^\sigma g = \sum_{n = 1}^\infty \lambda_n^\sigma P_n g, \ \
        \mathcal{D}(A^\sigma) = \left\{ g \in L^2(\Omega) \ \bigg| \ \sum_{n = 1}^\infty \lambda_n^{2\sigma} ||P_n g||_{L^2(\Omega)}^2 < \infty \right\}.
    \end{equation*}
    Then $||g||_{H^{2\sigma}(\Omega)} \leq C ||A^\sigma g||_{L^2(\Omega)}$ for $g \in \mathcal{D}(A^\sigma)$. Moreover, we denote $g(t) = O(t^{-\infty})$ if $g(t) = O(t^{-p})$ for any $p \in \mathbb{N}$.
    
    We now present the main theorems. Theorems~\ref{thm: IP1} and \ref{thm: IP2} specifically address issues (IP1) and (IP2), respectively. Additionally, we state the uniqueness of the whole source term $F(x, t) = f(x) \mu(t)$ in Theorem~\ref{thm: IP3}, which is comparable to Theorem 1.2 in \cite{jk22}.
    
    \begin{thm}
        \label{thm: IP1}
        Assume that $\mu \in L^\infty(0, \infty)$ has an asymptotic expansion over $t^{-1}$:
        \begin{equation}
            \mu(t) \sim \sum_{j = 0}^\infty \mu_j t^{-j} \text{ as } t \rightarrow \infty.
            \label{eq: expansion of mu}
        \end{equation}
        Besides, $\mu$ satisfies one of the following two conditions.
        \begin{itemize}
            \item Constant sign near infinity: 
            \begin{equation} \tag{4a}
                \inf \big\{ t \in [0, \infty] \ | \ \mu(s) \geq 0, \text{ a.e.} \, s \geq t \text{ or } \mu(s) \leq 0, \text{ a.e.} \, s \geq t \big\} < \infty.
                \label{eq: t_*}
            \end{equation}
            
            \item Fast decay: there exists $c_1$, $c_2 > 0$ such that
            \begin{equation} \tag{4b}
                |\mu(t)| \leq c_1 \exp(-c_2 \sqrt{t}), \ \forall t > 0.
                \label{eq: fast decay}
            \end{equation}
        \end{itemize}
        \addtocounter{equation}{1}
        (a) Assume that $f \in L^2(\Omega)$ and
        \begin{equation*}
            ||u(\cdot, t)||_{L^2(\omega)} = O(t^{-\infty}) \text{ as } t \rightarrow \infty.
        \end{equation*}
        Then $\mu = 0$ in $(0, \infty)$ or $f = 0$ in $\Omega$. \\
        (b) Let $f \in \mathcal{D}(A^{\sigma_2})$, where $\sigma_2 > (d - 1) / 4$. Assume that 
        \begin{equation*}
            ||\partial_\nu u(\cdot, t)||_{L^2(\gamma)} = O(t^{-\infty}) \text{ as } t \rightarrow \infty.
        \end{equation*}
        Then $\mu = 0$ in $(0, \infty)$ or $f = 0$ in $\Omega$.
    \end{thm}

    \begin{remark}
        \label{rem: constant sign}
        The assumption of constant sign is automatically satisfied when there exists $j \geq 0$ such that $\mu_j \neq 0$. Thus, the assumptions forces $\mu$ either to be non-oscillatory or to decay sufficiently fast when it has super-polynomial decay.
    \end{remark}
    
    \begin{thm}
        \label{thm: IP2}
        Let $\mu \in L^\infty(0, \infty)$ satisfy the conditions \eqref{eq: expansion of mu} and \eqref{eq: t_*}/\eqref{eq: fast decay} in Theorem~\ref{thm: IP1}. Define 
        \begin{equation}
            \Lambda(f) = \{ x_0 \in \Omega \ | \ P_n f(x_0) = 0 \text{ for all } n \in \mathbb{N} \}
            \label{eq: Lambda(f)}
        \end{equation}
        and
        \begin{equation}
            \Lambda_b(f) = \{ x_0 \in \partial \Omega \ | \ \partial_\nu P_n f(x_0) = 0 \text{ for all } n \in \mathbb{N} \}.
            \label{eq: Lambda_b(f)}
        \end{equation}
        (a) Let $f \neq 0$ and $f \in \mathcal{D}(A^{\sigma_3})$, where $\sigma_3 > \max \{ 0, d / 2 - 1 \}$. Then for any $x_0 \in \Omega \setminus \Lambda(f)$, 
        \begin{equation*}
            u(x_0, t) = O(t^{-\infty}) \text{ as } t \rightarrow \infty
        \end{equation*}
        implies $\mu = 0$ in $(0, \infty)$. \\
        (b) Let $f \neq 0$ and $f \in \mathcal{D}(A^{\sigma_4})$, where $\sigma_4 > d / 2 - 1 / 4$. Then for any $x_0 \in \partial \Omega \setminus \Lambda_b(f)$, 
        \begin{equation*}
            \partial_\nu u(x_0, t) = O(t^{-\infty}) \text{ as } t \rightarrow \infty
        \end{equation*}
        implies $\mu = 0$ in $(0, \infty)$.
    \end{thm}
    
    \begin{remark}
        \label{rem: Lambda(f)}
        For any nonzero function $f \in L^2(\Omega)$, $\Omega \setminus \Lambda(f)$ and $\partial \Omega \setminus \Lambda_b(f)$ are dense in $\Omega$ and $\partial \Omega$ respectively. This is due to the fact that $P_n f$ satisfies $(A - \lambda_n) P_n f = 0$ for any $n$. The unique continuation theorem for elliptic operators shows that if $P_n f$ or $\partial_\nu P_n f$ vanishes in an open set, then $P_n f = 0$ in $\Omega$. Therefore, if $\Lambda(f)$ or $\Lambda_b(f)$ has an interior point, then $f = 0$, which is a contradiction.
    \end{remark}
    
    \begin{thm}
        \label{thm: IP3}
        Let $\alpha \in (0, 2) \setminus \mathbb{Q}$. Suppose $u$ and $\widetilde{u}$ are solutions to the fractional diffusion-wave equation \eqref{eq: TFDE} with source terms $F(x, t) = f(x) \mu(t)$ and $\widetilde{F}(x, t) = \widetilde{f}(x) \widetilde{\mu}(t)$ respectively. Moreover, $\mu$ and $\widetilde{\mu} \in L^\infty(0, \infty)$ satisfy the conditions \eqref{eq: expansion of mu} and \eqref{eq: t_*}/\eqref{eq: fast decay} in Theorem~\ref{thm: IP1}. \\
        (a) Assume that $f, \widetilde{f} \in L^2(\Omega)$ and
        \begin{equation*}
            ||u(\cdot, t) - \widetilde{u}(\cdot, t)||_{L^2(\omega)} = O(t^{-\infty}) \text{ as } t \rightarrow \infty.
        \end{equation*}
        Then $F = \widetilde{F}$ in $\Omega \times (0, \infty)$. \\
        (b) Let $f, \widetilde{f} \in \mathcal{D}(A^{\sigma_2})$, where $\sigma_2 > (d - 1) / 4$. Assume that
        \begin{equation*}
            ||\partial_\nu u(\cdot, t) - \partial_\nu \widetilde{u}(\cdot, t)||_{L^2(\gamma)} = O(t^{-\infty}) \text{ as } t \rightarrow \infty.
        \end{equation*}
        Then $F = \widetilde{F}$ in $\Omega \times (0, \infty)$.
    \end{thm}

    \begin{remark}
        Our results cover the case when $\mu$ is compactly supported \cite{jk22,kly22,y23}. In this case $\mu$ has an asymptotic expansion of $0$, and both assumptions \eqref{eq: t_*} and \eqref{eq: fast decay} are satisfied. However, based on Laplace transform techniques in \cite{jk22,kly22}, the authors proved the uniqueness with lower regularity conditions, namely, $f \in L^2(\Omega)$ for analogs of Theorems \ref{thm: IP1} and \ref{thm: IP3}, and $f \in \mathcal{D} (A^{\sigma})$, $\sigma > (d - 2) / 4$ for an analog of Theorem \ref{thm: IP2}(b).
    \end{remark}
    
    \section{Proof of the main theorems}
    \label{sec: proof}
    
    \subsection{The key propositions}
    The proof of the main theorems rely on two key propositions. A critical component in this pursuit is the  integral function $\psi_n$, which plays a significant role in our analysis. We define it and establish its properties as follows:    

    \begin{lemma}
        \label{lem: psi_n}
        Define 
        \begin{equation}
            \psi_n(t) = \int_0^t (t - s)^{\alpha - 1} E_{\alpha, \alpha} (-\lambda_n (t - s)^\alpha) \mu(s) \dif{s}.
            \label{eq: psi_n}
        \end{equation}
        Then $\psi_n \in C[0, \infty)$ and the following estimate holds:
        \begin{equation*}
            |\psi_n(t)| \leq C(t) \, ||\mu||_{L^\infty (0, t)} \, \lambda_n^{-1} \log(1 + \lambda_n).
        \end{equation*}
    \end{lemma}

    \begin{pf}
        The estimate in Proposition~\ref{prop: Mittag-Leffler}(c) implies that $F(s) := s^{\alpha - 1} E_{\alpha, \alpha} (-\lambda_n s^\alpha) \in L^1_{loc} [0, \infty)$. Thus, for any $t_1, t_2 \geq 0$ we have
        \begin{align*}
            |\psi_n(t_1) - \psi_n(t_2)| 
            & \leq \int_0^{t_1} |F(t_2 - s) - F(t_1 - s)| \, |\mu(s)| \dif{s} + \int_{t_1}^{t_2} |F(t_2 - s)| \, |\mu(s)| \dif{s} \\
            & \leq ||F(\cdot + (t_2 - t_1)) - F||_{L^1(0, t_1)} ||\mu||_{L^\infty(0, t_2)} + ||F||_{L^1(0, t_2 - t_1)} ||\mu||_{L^\infty(0, t_2)}
        \end{align*}
        The continuity of $\psi_n$ comes from the continuity of translation under $L^1$ norm and absolute continuity of integrals. Moreover, 
        \begin{align*}
            |\psi_n(t)| &< C \int_0^t \frac{(t - s)^{\alpha - 1}} {1 + \lambda_n (t - s)^\alpha} \dif{s} \cdot ||\mu||_{L^\infty (0, t)} \\
            &= C ||\mu||_{L^\infty (0, t)} \, \lambda_n^{-1} \log( 1 + \lambda_n t^\alpha) \\
            &= C ||\mu||_{L^\infty (0, t)} \, \lambda_n^{-1} \log(1 + \lambda_n) \left( 1 + \frac{\log( (1 + \lambda_n)^{-1} + \lambda_n (1 + \lambda_n)^{-1} t^\alpha)} {\log(1 + \lambda_n)} \right) \\
            &\leq C ||\mu||_{L^\infty (0, t)} \, \lambda_n^{-1} \log(1 + \lambda_n) \left( 1 + \frac{\log( (1 + \lambda_1)^{-1} + t^\alpha)} {\log(1 + \lambda_1)} \right).
        \end{align*}
        The last inequality follows from $\lambda_n \geq \lambda_1$. The proof is complete by treating the last term as a constant dependent on $t$.
    \end{pf}
    
    Now we can state the key propositions.
    
    \begin{prop}
        \label{prop: key prop 12}
        Suppose $\mu \not\equiv 0$ satisfies the conditions \eqref{eq: expansion of mu} and \eqref{eq: t_*}/\eqref{eq: fast decay} in Theorem~\ref{thm: IP1}. We assume that $\{ a_n \}_{n \in \mathbb{N}} \subseteq \real$ satisfies $\sum_{n = 0}^\infty |a_n \lambda_n^{-1} \log(1 + \lambda_n)| < \infty$. If
        \begin{equation}
            \sum_{n = 0}^\infty a_n \psi_n(t) = O(t^{-\infty}) \text{ as } t \rightarrow \infty,
            \label{eq: decay condition 1}
        \end{equation}
        then $a_n = 0$ for all $n \in \mathbb{N}$.
    \end{prop}
    
    In view of Lemma~\ref{lem: psi_n} and the assumption of $a_n$ in Proposition~\ref{prop: key prop 12}, the sum $\sum_{n = 0}^\infty a_n \psi_n(t)$ is well defined. 
    
    The idea of proving Proposition~\ref{prop: key prop 12} is similar to that of Proposition 2 in \cite{y23}. The main task is to calculate the asymptotic expansion of $\psi_n(t)$. Then we can get the expansion of $\sum_{n = 0}^\infty a_n \psi_n(t)$. By the decay property in the proposition, all coefficients should be zero. The final step is to apply Lemma~\ref{lem: a_n=0} stated below to conclude that $a_n = 0$.
    
    \begin{lemma}
        \label{lem: a_n=0}
        Suppose $\sum_{n = 0}^\infty |a_n \lambda_n^{-1} \log(1 + \lambda_n)| < \infty$, and there is a sequence $\{ l_k \}_{k \in \mathbb{N}} \in \real$ increasing to infinity. If $\sum_{n = 1}^\infty a_n \lambda_n^{-l_k} = 0$ for all $k \in \mathbb{N}$, then $a_n = 0$ for all $n \in \mathbb{N}$. 
    \end{lemma}

    This lemma is almost the same as Lemma 3 in \cite{y23}, and it can be proved in exactly the same way, so we omit the proof here.

    Proposition~\ref{prop: key prop 3} is a slight modification of Proposition~\ref{prop: key prop 12} so that we can apply it to Theorem~\ref{thm: IP3}. To this end, we denote the expansion coefficients of $\widetilde{\mu}$ by $\widetilde{\mu}_j$ and denote $\widetilde{\psi}_n (t)$ by replacing $\mu$ with $\widetilde{\mu}$ in equation $\eqref{eq: psi_n}$.
    
    \begin{prop}
        \label{prop: key prop 3}
        Let $\alpha \in (0, 2) \setminus \mathbb{Q}$. Suppose $\mu$. $\widetilde{\mu}$ are nonzero and satisfy the conditions \eqref{eq: expansion of mu} and \eqref{eq: t_*}/\eqref{eq: fast decay} in Theorem~\ref{thm: IP1}. We assume that $\{ a_n \}_{n \in \mathbb{N}}$ and $\{ \widetilde{a}_n \}_{n \in \mathbb{N}}$ are nonzero sequences satisfying $\sum_{n = 0}^\infty |a_n \lambda_n^{-1} \log(1 + \lambda_n)| < \infty$ and $\sum_{n = 0}^\infty |\widetilde{a}_n \lambda_n^{-1} \log(1 + \lambda_n)| < \infty$. If 
        \begin{equation}
            \sum_{n = 0}^\infty a_n \psi_n(t) - \sum_{n = 0}^\infty \widetilde{a}_n \widetilde{\psi}_n(t) = O(t^{-\infty}) \text{ as } t \rightarrow \infty,
            \label{eq: decay condition 2}
        \end{equation}
        then there exists a constant $\kappa$ such that $a_n = \kappa \widetilde{a}_n$ for all $n \in \mathbb{N}$ and $\widetilde{\mu} = \kappa \mu$.
    \end{prop}

    \subsection{The asymptotic expansion of $\psi_n(t)$}
    Before showing the expansion of $\psi_n(t)$, we introduce some notations. Hereinafter, we let $b$ (with various subscripts) represent a non-essential coefficient that does not depend on $\lambda_n$ or $t$. The notation $b^\times$ (with various subscripts) represents a coefficient that vanishes if $\mu_j = 0$ for all $j \geq 1$. The specific meaning of symbols $b$ and $b^\times$ will vary according to context. 

    Given that $\alpha \neq 1$, the set $\{ k \in \mathbb{N} \ | -1 + k \alpha \notin \mathbb{N} \cup \{0\} \}$ is infinite. We relabel this set in ascending order as $\{ l_k \}_{k \in \mathbb{N}}$.
    For an arbitrary natural number $m$, we present an alternative expansion of $\mu$:
    \begin{equation}
        \mu(s) \sim \sum_{j = 1}^m \mu_j s^{-j} + \mu_{m + 1} s^{-m} (1 + s)^{-1} + \sum_{j = m + 2}^\infty (\mu_j + (-1)^{j - m} \mu_{m + 1}) s^{-j}.
        \label{eq: alternative expansion of mu}
    \end{equation}    
    Define $\overline{R}_{\mu, m'}(s)$ for $m' > m$ as
    \begin{equation}
        \overline{R}_{\mu, m'}(s) = \mu(s) - \sum_{j = 1}^m \mu_j s^{-j} - \mu_{m + 1} s^{-m} (1 + s)^{-1} - \sum_{j = m + 2}^{m'} (\mu_j + (-1)^{j - m} \mu_{m + 1}) s^{-j}.
        \label{eq: R_bar}
    \end{equation}
    Since $\overline{R}_{\mu, m + 1}(s) = O(s^{-m})$ as $s \rightarrow 0$ and $\overline{R}_{\mu, m + 1}(s) = O(s^{-m-2})$ as $s \rightarrow \infty$, the quantity 
    \begin{equation}
        \label{eq: c_mu, m} 
        c_{\mu, m} := \int_0^\infty s^m \overline{R}_{\mu, m + 1}(s) \dif{s}
    \end{equation}
    is well defined.

    \begin{lemma}
        \label{lem: expansion of psi_n}
        In the framework of Proposition~\ref{prop: key prop 12}, for any $N \in \mathbb{N}$, one can find $J$, $K$ and $M$ depending on $\alpha$ and $N$ only, such that $\psi_n(t)$ expands as follows:
        \begin{equation}
        \begin{aligned}
            \psi_n(t) &= \mu_0 \lambda_n^{-1} + \sum_{k = 1}^K \frac{(-1)^{l_k} \mu_0}{\Gamma(1 - l_k \alpha)} \lambda_n^{-l_k - 1} t^{-l_k \alpha} \\
            & \quad \ + \sum_{j = 1}^J \sum_{k = 1}^K \frac{(-1)^{l_k - 1 + j}} {\Gamma(-l_k \alpha)} \binom{- l_k \alpha - 1}{j - 1} \lambda_n^{-l_k - 1} t^{-l_k \alpha - j} (\mu_j \log t + c_{\mu, j - 1} + b_{jk}^\times) \\
            & \quad \ + \sum_{j = 1}^J \sum_{m = 0}^M b_{jm}^\times \lambda_n^{-1 - m / \alpha} t^{-j - m}
            + \lambda_n^{-1} \log(1 + \lambda_n) \, O(t^{-N}), \text{ as } t \rightarrow \infty.
            \label{eq: psi_n final}
        \end{aligned}
        \end{equation}
    \end{lemma}

    \noindent We emphasize that the bounding constant in the big-O notation is independent of $\lambda_n$ throughout the article.
    
    \begin{pf}
        {\bf Step 1:} Treatment of the zeroth order term $\mu_0$. 
        
        Since $\psi_n$ is linear with respect to $\mu$, we first consider the zeroth order term $\mu_0$. Note that $E_{\alpha, 1}(z)' = \alpha^{-1} E_{\alpha, \alpha}(z)$ by definition, so
        \begin{equation*}
            \int_0^t (t - s)^{\alpha - 1} E_{\alpha, \alpha} (-\lambda_n (t - s)^\alpha) \mu_0 \dif{s} 
                            = \mu_0 \lambda_n^{-1} (1 - E_{\alpha, 1}(-\lambda_n t^\alpha)).
        \end{equation*}
        Recall the definition of $l_k$, and the asymptotic expansion of the Mittag-Leffler function in Proposition~\ref{prop: Mittag-Leffler}(b) becomes
        \begin{equation*}
            E_{\alpha, 1}(-\lambda_n t^\alpha) \sim \sum_{k = 1}^\infty \frac{(-1)^{l_k + 1}} {\Gamma(1 - l_k \alpha)} (\lambda_n t^\alpha)^{-l_k}.
        \end{equation*}
        Fix $N \in \mathbb{N}$, then we can choose a sufficiently large $K = K(N)$ such that $\alpha l_{K + 1} > N$, and we have
        \begin{equation*}
                \int_0^t (t - s)^{\alpha - 1} E_{\alpha, \alpha} (-\lambda_n (t - s)^\alpha) \mu_0 \dif{s} 
                = \mu_0 \lambda_n^{-1} + \sum_{k = 1}^K \frac{(-1)^{l_k} \mu_0}{\Gamma(1 - l_k \alpha)} \lambda_n^{-l_k - 1} t^{-l_k \alpha} + \lambda_n^{-1} O(t^{-N}).
        \end{equation*}

        \noindent{\bf Step 2:} Expanding $\mu$ away from 0.
        
        In the remaining proof, we set $\mu_0 = 0$ and introduce two parameters $\beta_1$ and $\beta_2$ so that $0 < \beta_1 < \beta_2 < 1$. Write $\mu(t) = \sum_{j = 1}^J \mu_j t^{-j} + R_{\mu, J}(t)$, then $|R_{\mu, J}(t)| < Ct^{-J - 1}$ for $t > 1$. Since
        \begin{equation*}
            \psi_n(t) = \int_0^{t^{\beta_1}} (t - s)^{\alpha - 1} E_{\alpha, \alpha} (-\lambda_n (t - s)^\alpha) \mu(s) \dif{s} + \int_0^{t - t^{\beta_1}} s^{\alpha - 1} E_{\alpha, \alpha} (-\lambda_n s^\alpha) \mu(t - s) \dif{s} 
        \end{equation*}
        and
        \begin{align*}
            & \quad \ \left| \int_0^{t - t^{\beta_1}} s^{\alpha - 1} E_{\alpha, \alpha} (-\lambda_n s^\alpha) \, R_{\mu, J}(t - s) \dif{s} \right| \\
            & \leq \int_0^{t - t^{\beta_1}} s^{\alpha - 1} \frac{C}{1 + \lambda_n s^\alpha} \cdot C(t - s)^{-J - 1} \dif{s} \\
            & \leq C \lambda_n^{-1} \log(1 + \lambda_n t^\alpha) t^{-(J + 1) \beta_1} \\
            & \leq C \lambda_n^{-1} \log(1 + \lambda_n) (1 + C \log t) t^{-(J + 1) \beta_1}
        \end{align*}
        for $t > 1$, choosing a sufficiently large $J = J(N)$ such that $(J + 1) \beta_1 > N$, we get
        \begin{equation}
        \begin{aligned}
            \psi_n(t) = &\int_0^{t^{\beta_1}} (t - s)^{\alpha - 1} E_{\alpha, \alpha} (-\lambda_n (t - s)^\alpha) \mu(s) \dif{s} \\
            & + \sum_{j = 1}^J \mu_j \int_0^{t - t^{\beta_1}} s^{\alpha - 1} E_{\alpha, \alpha} (-\lambda_n s^\alpha) (t - s)^{-j} \dif{s}
            + \lambda_n^{-1} \log(1 + \lambda_n) \, O(t^{-N}).
            \label{eq: psi_n 1} 
        \end{aligned}
        \end{equation}

        \noindent{\bf Step 3:} Expanding $E_{\alpha, \alpha} (-\lambda_n s^\alpha)$ for $s$ away from 0.
        
        By definition of $l_k$, we can find that the set $\{ k \in \mathbb{N} \ | -\alpha + k \alpha \notin \mathbb{N} \cup \{ 0 \} \}$ coincides with $\{ l_k + 1 \}_{k \in \mathbb{N}}$. Hence, for $t > 1$ and $s > \lambda_n^{-1 / \alpha} t^{\beta_2}$, we have
        \begin{equation}
            E_{\alpha, \alpha} (-\lambda_n s^\alpha) = \sum_{k = 1}^K \frac{(-1)^{l_k}}{\Gamma(-l_k \alpha)} (\lambda_n s^\alpha)^{-l_k - 1} + R_{E, K}(\lambda_n s^\alpha)
            \label{eq: MLE expansion}
        \end{equation}
        with $R_{E, K}(\lambda_n s^\alpha) = O\big( (\lambda_n s^\alpha)^{-l_{K + 1} - 1} \big)$. We can choose $T_1$ sufficiently large and independent of $\lambda_n$ such that $t - t^{\beta_1} > \lambda_n^{-1 / \alpha} t^{\beta_2}$ when $t > T_1$. This guarantees that, when $t > T_1$ and $0 \leq s \leq t^{\beta_1}$, we have $\lambda_n (t - s)^\alpha > T_1^{\alpha \beta_2}$, which is independent of $n$. Together with $\mu \in L^\infty (0, \infty)$, we get
        \begin{equation}
        \begin{aligned}
            &\quad \ \left| \int_0^{t^{\beta_1}} (t - s)^{\alpha - 1} R_{E, K}(\lambda_n (t - s)^\alpha) \mu(s) \dif{s} \right| \\
            &\leq \int_0^{t^{\beta_1}} (t - s)^{\alpha - 1} \cdot C (\lambda_n (t - s)^\alpha)^{-l_{K + 1} - 1} |\mu(s)| \dif{s} \\
            &= \lambda_n^{-l_{K + 1} - 1} ||\mu||_{L^\infty(0, \infty)} \, O(t^{-l_{K + 1} \alpha}).
            \label{eq: R_E_K 1}
        \end{aligned}
        \end{equation}
        Besides, we have $\lambda_n s^\alpha > T_1^{\alpha \beta_2}$ for $t > T_1$ and $\lambda_n^{-1 / \alpha} t^{\beta_2} \leq s \leq t - t^{\beta_1}$. Hence, for each $j \in \mathbb{N}$, 
        \begin{equation}
        \begin{aligned}
            &\quad \ \left| \int_{\lambda_n^{-1 / \alpha} t^{\beta_2}}^{t - t^{\beta_1}} s^{\alpha - 1} R_{E, K}(\lambda_n s^\alpha) (t - s)^{-j} \dif{s} \right| \\
            &\leq \int_{\lambda_n^{-1 / \alpha} t^{\beta_2}}^{t - t^{\beta_1}} s^{\alpha - 1} \cdot C (\lambda_n s^\alpha)^{-l_{K + 1} - 1} (t - s)^{-j} \dif{s} \\
            &\leq C \lambda_n^{-l_{K + 1} - 1} t^{-\beta_1 j} \int_{\lambda_n^{-1 / \alpha} t^{\beta_2}}^\infty s^{-1 - l_{K + 1} \alpha} \dif{s} \\
            &= \lambda_n^{-1} O(t^{-\beta_1 j - \alpha \beta_2 l_{K + 1} }).
            \label{eq: R_E_K 2}
        \end{aligned}
        \end{equation}
        For $N \in \mathbb{N}$ given above, we choose $K = K(N)$ such that $\alpha \beta_2 l_{K + 1} > N$. Combining equations \eqref{eq: psi_n 1} -- \eqref{eq: R_E_K 2}, we obtain
        \begin{equation}
        \begin{aligned}
            \psi_n(t) = &\sum_{k = 1}^K \frac{(-1)^{l_k}}{\Gamma(-l_k \alpha)} \lambda_n^{-l_k - 1} \int_0^{t^{\beta_1}} (t - s)^{-l_k \alpha - 1} \mu(s) \dif{s} \\
            & + \sum_{j = 1}^J \sum_{k = 1}^K \mu_j \frac{(-1)^{l_k}}{\Gamma(-l_k \alpha)} \lambda_n^{-l_k - 1} \int_{ \lambda_n^{-1 / \alpha} t^{\beta_2} }^{t - t^{\beta_1}} s^{- 1 - l_k \alpha} (t - s)^{-j} \dif{s} \\
            & + \sum_{j = 1}^J \mu_j \int_0^{\lambda_n^{-1 / \alpha} t^{\beta_2}} s^{\alpha - 1} E_{\alpha, \alpha} (-\lambda_n s^\alpha) (t - s)^{-j} \dif{s} \\
            & + \lambda_n^{-1} \log(1 + \lambda_n) \, O(t^{-N}).
            \label{eq: psi_n 2}
        \end{aligned}
        \end{equation}

        \noindent{\bf Step 4:} Expanding $(t - s)^\gamma, \gamma < 0$ in the integrals.
        
        Since
        \begin{equation*}
            (1 - s)^\gamma = \sum_{m = 0}^M (-1)^m \binom{\gamma}{m} s^m + R_{\gamma, M}(s),
        \end{equation*}
        we obtain that
        \begin{equation*}
            (t - s)^\gamma = t^\gamma \sum_{m = 0}^M (-1)^m \binom{\gamma}{m} s^m t^{-m} + R_{\gamma, M}(s/t). 
        \end{equation*}
        We can choose $T_2 > 0$ independent of $\lambda_n$ such that $t^{\beta_1 - 1} < 1/2$ and $\lambda_n^{-1 / \alpha} t^{\beta_2 - 1} < 1/2$ for $t > T_2$. When $t > T_2$ and $0 \leq s \leq \max \{ t^{\beta_1}, \lambda_n^{-1 / \alpha} t^{\beta_2} \}$, we have $R_{\gamma, M}(s/t) = O( (s/t)^{M + 1} )$. It follows that
        \begin{equation}
        \begin{aligned}
            &\quad \ \int_0^{t^{\beta_1}} (t - s)^{-l_k \alpha - 1} \mu(s) \dif{s} \\
            &= \sum_{m = 0}^M (-1)^m \binom{-l_k \alpha - 1}{m} t^{-l_k \alpha - 1 - m} \int_0^{t^{\beta_1}} s^m \mu(s) \dif{s} + t^{-l_k \alpha - 1} \int_0^{t^{\beta_1}} R_{-l_k \alpha - 1, M}(s/t) \mu(s) \dif{s} \\
            &= \sum_{m = 0}^M (-1)^m \binom{-l_k \alpha - 1}{m} t^{-l_k \alpha - 1 - m} \int_0^{t^{\beta_1}} s^m \mu(s) \dif{s} + O(t^{-l_k \alpha - 1 - (M + 1)(1 - \beta_1)}).
            \label{eq: R_gamma_M 1}
        \end{aligned}
        \end{equation}
        On the other hand, 
        \begin{equation*}
        \begin{aligned}
            &\quad \ \int_0^{\lambda_n^{-1 / \alpha} t^{\beta_2}} s^{\alpha - 1} E_{\alpha, \alpha}(-\lambda_n s^\alpha) R_{-j, M}(s/t) \dif{s} \\
            &\leq C \int_0^{\lambda_n^{-1 / \alpha} t^{\beta_2}} s^{\alpha - 1} (1 + \lambda_n s^\alpha)^{-1} (s / t)^{M + 1} \dif{s} \\
            &\leq C \lambda_n^{-1 - (M + 1) / \alpha} \, t^{- M - 1} \int_0^{t^{\alpha \beta_2}} (1 + s)^{-1} s^{(M + 1) / \alpha} \dif{s} \\
            &\leq C \lambda_n^{-1 - (M + 1) / \alpha} \, t^{- M - 1} \int_0^{t^{\alpha \beta_2}} s^{(M + 1) / \alpha - 1} \dif{s} \\
            &= \lambda_n^{-1 - (M + 1) / \alpha} \, O( t^{-(M + 1)(1 - \beta_2)} ).
        \end{aligned}
        \end{equation*}
        It follows that 
        \begin{equation}
        \begin{aligned}
            &\quad \ \int_0^{\lambda_n^{-1 / \alpha} t^{\beta_2}} s^{\alpha - 1} E_{\alpha, \alpha} (-\lambda_n s^\alpha) (t - s)^{-j} \dif{s} \\
            &= \sum_{m = 0}^M (-1)^m \binom{-j}{m} t^{-j - m} \int_0^{\lambda_n^{-1 / \alpha} t^{\beta_2}} s^{\alpha - 1 + m} E_{\alpha, \alpha} (-\lambda_n s^\alpha) \dif{s} \\
            &\qquad \quad + t^{-j} \int_0^{\lambda_n^{-1 / \alpha} t^{\beta_2}} s^{\alpha - 1} E_{\alpha, \alpha} (-\lambda_n s^\alpha) R_{-j, M}(s/t) \dif{s} \\
            &= \sum_{m = 0}^M \frac{(-1)^m}{\alpha} \binom{-j}{m} \lambda_n^{-1 - m / \alpha} t^{-j - m} \int_0^{t^{\alpha \beta_2}} s^{m / \alpha} E_{\alpha, \alpha} (-s) \dif{s} \\
            &\qquad \quad + \lambda_n^{-1 - (M + 1) / \alpha} \, O( t^{-(M + 1)(1 - \beta_2)} ). \\
            \label{eq: R_gamma_M 2}
        \end{aligned}
        \end{equation}
        Choose $M = M(N)$ such that $1 + (M + 1)(1 - \beta_1) > N$ and $(M + 1)(1 - \beta_2) > N$. Then by equations \eqref{eq: psi_n 2}, \eqref{eq: R_gamma_M 1} and \eqref{eq: R_gamma_M 2}, we get
        \begin{equation}
        \begin{aligned}
            \psi_n(t) = &\sum_{k = 1}^K \sum_{m = 0}^M \frac{(-1)^{l_k + m}}{\Gamma(-l_k \alpha)} \binom{-l_k \alpha - 1}{m} \lambda_n^{-l_k - 1} t^{-l_k \alpha - 1 - m} \int_0^{t^{\beta_1}} s^m \mu(s) \dif{s} \\
            & + \sum_{j = 1}^J \sum_{k = 1}^K \mu_j \frac{(-1)^{l_k}}{\Gamma(-l_k \alpha)} \lambda_n^{-l_k - 1} t^{-l_k \alpha - j} \int_{ \lambda_n^{-1 / \alpha} t^{\beta_2 - 1} }^{1 - t^{\beta_1 - 1}} s^{-l_k \alpha - 1} (1 - s)^{-j} \dif{s} \\
            & + \sum_{j = 1}^J \sum_{m = 0}^M \frac{(-1)^m \mu_j}{\alpha} \binom{-j}{m} \lambda_n^{-1 - m / \alpha} t^{-j - m} \int_0^{t^{\alpha \beta_2}} s^{m / \alpha} E_{\alpha, \alpha} (-s) \dif{s} \\
            & + \lambda_n^{-1} \log(1 + \lambda_n) \, O(t^{-N}).
            \label{eq: psi_n 3}
        \end{aligned}
        \end{equation}

        \noindent{\bf Step 5:} Expanding integrals in the form of $\int s^\gamma F(s) \dif{s}$.

        \noindent{\bf Step 5.1:} Expansion of $\int_0^{t^{\beta_1}} s^m \mu(s) \dif{s}$.
        
        We use the alternative expansion \eqref{eq: alternative expansion of mu} of $\mu$ to expand the integral $\int_0^{t^{\beta_1}} s^m \mu(s) \dif{s}$. By equation \eqref{eq: R_bar}, $\overline{R}_{\mu, J'} = O(t^{-J' - 1})$ for $J' > m$. Hence, 
        \begin{align*}
            \int_0^{t^{\beta_1}} s^m \mu(s) \dif{s}
            &= \sum_{j = 1}^m \int_0^{t^{\beta_1}} \mu_j s^{m - j} \dif{s} + \int_0^{t^{\beta_1}} \mu_{m + 1} (1 + s)^{-1} \dif{s} + \int_0^\infty s^m \overline{R}_{\mu, m + 1}(s) \dif{s} \\
            &\quad \ - \sum_{j = m + 2}^{J'} \int_{t^{\beta_1}}^\infty (\mu_j + (-1)^{j - m} \mu_{m + 1}) s^{m - j} \dif{s} - \int_{t^{\beta_1}}^\infty s^m \overline{R}_{\mu, J'}(s) \dif{s} \\
            &= \sum_{j = 1}^m \frac{\mu_j}{m - j + 1} t^{(m - j + 1) \beta_1} + \mu_{m + 1} \log(1 + t^{\beta_1}) + c_{\mu, m} \\
            &\quad \ - \sum_{j = m + 2}^{J'} \frac{\mu_j + (-1)^{j - m} \mu_{m + 1}}{m - j + 1} t^{(m - j + 1) \beta_1} + O(t^{(m - J') \beta_1}).
        \end{align*}
        Choosing sufficiently large $J' = J'(N)$ such that $(J' - M) \beta_1 > N$, we obtain
        \begin{equation}
        \begin{aligned}
            \int_0^{t^{\beta_1}} s^m \mu(s) \dif{s} 
            &= \sum_{\substack{1 \leq j \leq J' \\ j \neq m + 1}} b_j^\times t^{(m - j + 1) \beta_1} + \mu_{m + 1} \log(1 + t^{\beta_1}) + c_{\mu, m} + O(t^{-N}) \\
            &= \sum_{\substack{1 \leq j \leq J' \\ j \neq m + 1}} b_j^\times t^{(m - j + 1) \beta_1} + \mu_{m + 1} \beta_1 \log t + c_{\mu, m} + O(t^{-N}).
            \label{eq: first term}
        \end{aligned}
        \end{equation}

        \noindent{\bf Step 5.2:} Expansion of $\int_0^{t^{\alpha \beta_2}} s^{m / \alpha} E_{\alpha, \alpha} (-s) \dif{s}$.
        
        We treat the integral $\int_0^{t^{\alpha \beta_2}} s^{m / \alpha} E_{\alpha, \alpha} (-s) \dif{s}$ similarly. When $m / \alpha \in \mathbb{N} \cup \{0\}$, observe that the coefficient of $s^{-m / \alpha - 1}$ in the expansion of $E_{\alpha, \alpha}(-s)$ is $0$, so for any $m \in \mathbb{N} \cup \{0\}$, we can use the ordinary expansion of $E_{\alpha, \alpha}$:
        \begin{equation*}
            E_{\alpha, \alpha}(-s) = \sum_{k = 1}^{K'} \frac{(-1)^{l_k}}{\Gamma(-l_k \alpha)} s^{-l_k - 1} + R_{E, K'}(s),
        \end{equation*}
        where $R_{E, K'}$ is defined in equation \eqref{eq: MLE expansion}. By a similar calculation as above we have
        \begin{equation}
            \int_0^{t^{\alpha \beta_2}} s^{m / \alpha} E_{\alpha, \alpha} (-s) \dif{s}
            = \sum_{k = 1}^{K'} b_{j, k} t^{(m - l_k \alpha) \beta_2} + b_{j, m} + O(t^{-N})
            \label{eq: third term}
        \end{equation}
        for some $K' = K'(N)$ satisfying $(l_{K + 1} \alpha - M) \beta_2 > N$.

        \noindent{\bf Step 5.3:} Expansion of $\int_{ \lambda_n^{-1 / \alpha} t^{\beta_2 - 1} }^{1 - t^{\beta_1 - 1}} s^{-l_k \alpha - 1} (1 - s)^{-j} \dif{s}$.
        
        This integral can be expressed by the incomplete Beta function $B(a, b; x)$ and the hypergeometric function ${}_2\mathbf{F}_1(a, b, c; x)$ (see Section 8.17 and Chapter 15 of \cite{handbook}). By equations (8.17.7) and (15.1.2) of \cite{handbook} we obtain
        \begin{equation}
        \begin{aligned}
            &\quad \ \int_{ \lambda_n^{-1 / \alpha} t^{\beta_2 - 1} }^{1 - t^{\beta_1 - 1}} s^{-l_k \alpha - 1} (1 - s)^{-j} \dif{s} \\
            &= B(-l_k \alpha, 1 - j; 1 - t^{\beta_1 - 1}) - B(-l_k \alpha, 1 - j; \lambda_n^{-1 / \alpha} t^{\beta_2 - 1}) \\
            &= \Gamma(-l_k \alpha) (1 - t^{\beta_1 - 1})^{-l_k \alpha} \ {}_2\mathbf{F}_1(-l_k \alpha, j, 1 - l_k \alpha; 1 - t^{\beta_1 - 1}) \\
            &\quad \ - \Gamma(-l_k \alpha) (\lambda_n^{-1 / \alpha} t^{\beta_2 - 1})^{-l_k \alpha} \ {}_2\mathbf{F}_1(-l_k \alpha, j, 1 - l_k \alpha; \lambda_n^{-1 / \alpha} t^{\beta_2 - 1}).
            \label{eq: second term 1}
        \end{aligned}
        \end{equation}
        The first hypergeometric function is expanded via equations (15.8.10) and (15.8.12) in \cite{handbook}. After simplifying the expression we get
        \begin{equation}
        \begin{aligned}
            &\quad \ {}_2\mathbf{F}_1 (-l_k \alpha, j, 1 - l_k \alpha; 1 - t^{\beta_1 - 1}) \\
            &= \sum_{m = 0}^\infty b_m t^{(\beta_1 - 1) (m - j + 1)} - \frac{(-1)^{j - 1}}{\Gamma(-l_k \alpha - j + 1)} \sum_{m = 0}^\infty \frac{(-l_k \alpha)_m (j)_m} {(m + j - 1)! \, m!} t^{(\beta_1 - 1)m} \log t^{\beta_1 - 1} \\
            &= \sum_{m = 0}^\infty b_m t^{(\beta_1 - 1) (m - j + 1)} + \frac{(-1)^{j - 1} (1 - \beta_1)} {\Gamma(-l_k \alpha - j + 1) (j - 1)!} (1 - t^{\beta_1 - 1})^{l_k \alpha} \log t.
            \label{eq: hypergeometric 1}
        \end{aligned}
        \end{equation}
        For the second hypergeometric function we expand it via its definition
        \begin{equation}
            \quad \ {}_2\mathbf{F}_1 (-l_k \alpha, j, 1 - l_k \alpha; \lambda_n^{-1 / \alpha} t^{\beta_2 - 1}) = \sum_{m = 0}^\infty b_m (\lambda_n^{-1 / \alpha} t^{\beta_2 - 1})^m.
            \label{eq: hypergeometric 2}
        \end{equation}
        Equations \eqref{eq: second term 1}, \eqref{eq: hypergeometric 1} and \eqref{eq: hypergeometric 2} imply
        \begin{align*}
            &\quad \ \int_{ \lambda_n^{-1 / \alpha} t^{\beta_2 - 1} }^{1 - t^{\beta_1 - 1}} s^{-l_k \alpha - 1} (1 - s)^{-j} \dif{s} \\
            &= (1 - t^{\beta_1 - 1})^{-l_k \alpha} \sum_{m = 0}^\infty b_m t^{(\beta_1 - 1) (m - j + 1)} + (-1)^{j - 1} \binom{-l_k \alpha - 1}{j - 1} (1 - \beta_1) \log t \\
            &\quad \ - \sum_{m = 0}^\infty b_m (\lambda_n^{-1 / \alpha} t^{\beta_2 - 1})^{-l_k \alpha + m} \\
            &= \sum_{m = 0}^\infty b_m t^{(\beta_1 - 1) (m - j + 1)} + (-1)^{j - 1} \binom{-l_k \alpha - 1}{j - 1} (1 - \beta_1) \log t - \sum_{m = 0}^\infty b_m (\lambda_n^{-1 / \alpha} t^{\beta_2 - 1})^{-l_k \alpha + m}.
        \end{align*}
        We can choose $M' = M'(N)$ such that $(1 - \beta_1)(M' + 2) > N$ and $(1 - \beta_2)(M' + 1) > N$, then we have
        \begin{equation}
        \begin{aligned}
            &\quad \ \int_{ \lambda_n^{-1 / \alpha} t^{\beta_2 - 1} }^{1 - t^{\beta_1 - 1}} s^{-l_k \alpha - 1} (1 - s)^{-j} \dif{s} \\
            &= \sum_{m = 0}^{M'} b_m t^{(\beta_1 - 1) (m - j + 1)} + (-1)^{j - 1} \binom{-l_k \alpha - 1}{j - 1} (1 - \beta_1) \log t \\
            &\quad \ - \sum_{m = 0}^{M'} b_m (\lambda_n^{-1 / \alpha} t^{\beta_2 - 1})^{-l_k \alpha + m} + O \big( t^{-N + j} + \lambda_n^{l_k} t^{-N + l_k \alpha} \big).
            \label{eq: second term 2}
        \end{aligned}
        \end{equation}
        
        Now we can substitute equations \eqref{eq: first term}, \eqref{eq: third term} and \eqref{eq: second term 2} into \eqref{eq: psi_n 3} and simplify the expression. Since adding more terms to the partial expansion will keep the remainder $\lambda_n^{-1} \log(1 + \lambda_n) \, O(t^{-N})$ unchanged, by denoting $J'' = \max \{ J, J', M+ 1 \}$, $K'' = \max \{ K, K' \}$ and $M'' = \max \{ M, M' \}$, we have
        \begin{align*}
            \psi_n(t) = &\sum_{j = 1}^{J''} \sum_{k = 1}^{K''} \sum_{m = 0}^{M''} \Big( b_{jkm}^\times \lambda_n^{-l_k - 1} t^{-l_k \alpha - \beta_1 j - (1 - \beta_1)(m + 1)} + b_{jkm}^\times \lambda_n^{-1 - m / \alpha} t^{-\beta_2 l_k \alpha - j - (1 - \beta_2) m} \Big) \\
            & + \sum_{j = 1}^{J''} \sum_{k = 1}^K \frac{(-1)^{l_k - 1 + j}} {\Gamma(-l_k \alpha)} \binom{-l_k \alpha - 1}{j - 1} \lambda_n^{-l_k - 1} t^{-l_k \alpha - j} (\mu_j \log t + c_{\mu, j - 1}) \\
            & + \sum_{j = 1}^J \sum_{m = 0}^M b_{jm}^\times \lambda_n^{-1 - m / \alpha} t^{-j - m}
            + \lambda_n^{-1} \log(1 + \lambda_n) \, O(t^{-N}).
        \end{align*}

        \noindent{\bf Step 6:} Independence of the expansion from $\beta_1$ and $\beta_2$.
        
        Finally we claim that the expansion of $\psi_n$ should be independent of $\beta_1$ and $\beta_2$ due to the uniqueness of asymptotic expansion. We choose $(\beta_1, \beta_2)$ and $(\beta_1', \beta_2')$ such that $\beta_1$, $\beta_2$, $\beta_1'$, $\beta_2'$ and $\alpha$ are algebraically independent over $\mathbb{Z}$. It follows that, if we pick a term in the first sum of the expansion of $\psi_n(t)$ with parameters $(\beta_1, \beta_2)$, the term still occurs under parameters $(\beta_1', \beta_2')$ if and only if the coefficients of $\beta_1$ and $\beta_2$ on the exponent of $t$ vanishes. Note that $l_k \alpha \notin \mathbb{N} \cup \{ 0 \}$ by definition, so the only remaining term in the first sum is $\sum_{j = 1}^{J''} \sum_{k = 1}^{K''} b_{jk, j-1}^\times \lambda_n^{-l_k - 1} t^{-l_k \alpha - j}$. This term can be brought into the second sum, again due to $l_k \alpha \notin \mathbb{N} \cup \{ 0 \}$.
    \end{pf}
    
    \subsection{Proof of the key propositions}
    \begin{lemma}
        \label{lem: mu equiv 0}
        Suppose $\mu \in L^\infty (0, \infty)$ satisfies the conditions \eqref{eq: expansion of mu} and \eqref{eq: t_*}/\eqref{eq: fast decay} in Theorem \ref{thm: IP1}. Define $c_{\mu, m}$ according to equations \eqref{eq: R_bar} and \eqref{eq: c_mu, m}. If $\mu_j = 0$ and $c_{\mu, m} = 0$ for all $j, m \geq 0$, then $\mu \equiv 0$.
    \end{lemma}

    \begin{proof}
        Since $\mu_j = 0$ for all $j \geq 0$, by definition \eqref{eq: R_bar} and \eqref{eq: c_mu, m} of $c_{\mu, m}$ we have $c_{\mu, m} = \int_0^\infty s^m \mu(s) \dif{s}$. 
        
        Suppose $\mu$ satisfies the constant sign condition \eqref{eq: t_*}, and we take $t_*$ to be the last change of sign of $\mu$ defined in equation \eqref{eq: t_*}. Without loss of generality, we assume that $\mu(t) \geq 0$ a.e. $\! t \geq t_*$. If $\mu \equiv 0$ on $t \geq t_*$, then we can further assume that $t_* > 0$. Now $c_{\mu, m} = 0$ implies $\int_0^{t_*} s^m \mu(s) \dif{s} = 0$. By Weierstrass approximation theorem, the set of polynomials are dense in $L^1(0, t_*)$. It follows that $\mu = 0$ on $(0, t_*)$. If $\mu \not\equiv 0$ on $t \geq t_*$, then we have
        \begin{equation*}
            \int_1^\infty s^m \mu(s / t_*) \dif{s} = -\int_0^1 s^m \mu(s / t_*) \dif{s}.
        \end{equation*}
        On the one hand, 
        \begin{equation*}
            \lim_{m \rightarrow \infty} \int_0^1 s^m \mu(s / t_*) \dif{s} = 0,
        \end{equation*}
        On the other hand, there exists $\delta > 0$ such that $\mu \not\equiv 0$ on $t > t_* + \delta$. Since $\mu \geq 0$ on $t > t_* + \delta$, we can bound $\mu |_{(t_* + \delta, \infty)}$ with a simple function $c_0 \chi_E$ from below, where $c_0 > 0$ and $E$ is a set of positive measure $|E|$. It follows that 
        \begin{equation*}
            \int_1^\infty s^m \mu(s / t_*) \dif{s} \geq \int_{E / t_*} c_0 s^m \dif{s} \geq c_0 (1 + \delta / t_*)^m |E|.
        \end{equation*}
        Thus, 
        \begin{equation*}
            \lim_{m \rightarrow \infty} \int_1^\infty s^m \mu(s / t_*) \dif{s} = \infty.
        \end{equation*}
        which is a contradiction.
        
        Now we consider the case where $\mu$ satisfies the decay condition \eqref{eq: fast decay}. We set $\mu^{(1)} (t) = c_1 \exp(- c_2 \sqrt{t} / 2)$ and $\mu^{(2)} = \mu / \mu^{(1)}$, then we have $|\mu^{(2)} (t)| \leq \exp(-c_2 \sqrt{t} / 2)$. We claim that the set $\{ p \mu^{(1)} \, |$ $p$ is a polynomial$\}$ is dense on $C_0[0, \infty)$. By Theorems~\ref{thm: IP1} and \ref{thm: IP2} of \cite{p55}, it is equivalent to prove that
        \begin{equation}
            \sup_p \int_{\real} \frac{ \log_+ |p(x)| }{1 + x^2} \dif{x} = \infty,
            \label{eq: density condition}
        \end{equation}
        where $\log_+ x = \max\{ \log x, 0 \}$ and the supremum is taken over all polynomials $p$ such that $|p(x)| \leq c_1 \exp(c_2 |x| / 2)$. We construct a (pointwise) increasing sequence of polynomials $p_n$ by
        \begin{equation*}
            p_n(x) = c_1 \sum_{k = 0}^n \frac{(\frac{1}{2} c_2 x)^{2k}} {(2k)!}.
        \end{equation*}
        Note that
        \begin{equation*}
            \lim_{n \rightarrow \infty} p_n(x) = \frac{1}{2} c_1 \big( \exp(c_2 x / 2) + \exp(-c_2 x / 2) \big),
        \end{equation*}
        which is between $c_1 / 2 \cdot \exp(c_2 |x| / 2)$ and $c_1 \exp(c_2 |x| / 2) $. Since $\log_+$ is an increasing function, by monotone convergence theorem we have
        \begin{align*}
            \lim_{n \rightarrow \infty} \int_{\real} \frac{ \log_+ |p_n(x)| }{1 + x^2} \dif{x}
            & = \int_{\real} \frac{ \log_+ \left| \lim\limits_{n \rightarrow \infty} p_n(x) \right| }{1 + x^2} \dif{x} \\
            & \geq \int_{\real} \frac{ \log_+ \left( c_1 / 2 \cdot \exp(c_2 |x| / 2) \right) }{1 + x^2} \dif{x} \\
            & = \int_{|x| > c_3} \frac{ \log \left( c_1 / 2 \right) + c_2 |x| / 2 }{1 + x^2} \dif{x} \\
            & = \infty
        \end{align*}
        So the condition \eqref{eq: density condition} is satisfied. Since $c_{\mu, m} = 0$ for all $m \geq 0$, the equation $\int_0^\infty p \mu^{(1)} \cdot \mu^{(2)} \dif{s} = 0$ is true for all polynomials $p$. The decay property $|\mu^{(2)} (t)| \leq \exp(-c_2 \sqrt{t} / 2)$ implies that $\mu^{(2)} \in L^1(0, \infty)$, which further implies $\mu^{(2)} \dif{s} \in \left( C_0[0, \infty) \right)^*$. Finally we can use the density of the set $\{ p \mu^{(1)} \, | \ p \text{ is a polynomial} \}$ to obtain $\mu^{(2)} \equiv 0$, which follows $\mu \equiv 0$.
    \end{proof}
    
    \begin{pf}[Proposition~\ref{prop: key prop 12}]
        We have shown that the series $\sum_{n = 1}^\infty a_n \psi_n(t)$ is absolutely convergent. Write $A_\gamma = \sum_{n = 1}^\infty a_n \lambda_n^{-\gamma - 1}$, $\gamma \in \real$, then by Lemma~\ref{lem: expansion of psi_n} and the condition 
        $$
        \sum_{n = 1}^\infty |a_n \lambda_n^{-1} \log (1 + \lambda_n)| < \infty,
        $$
        the expansion of $\sum_{n = 1}^\infty a_n \psi_n(t)$ is obtained as follows:
        \begin{equation}
        \begin{aligned}
            \sum_{n = 1}^\infty a_n \psi_n(t) &= \mu_0 A_0 + \sum_{k = 1}^K \frac{(-1)^{l_k} \mu_0}{\Gamma(1 - l_k \alpha)} A_{l_k} t^{-l_k \alpha} \\
            & \quad \ + \sum_{j = 1}^J \sum_{k = 1}^K \frac{(-1)^{l_k - 1 + j}} {\Gamma(-l_k \alpha)} \binom{-l_k \alpha - 1}{j - 1} A_{l_k} t^{-l_k \alpha - j} (\mu_j \log t + c_{\mu, j - 1} + b_{jk}^\times) \\
            & \quad \ + \sum_{j = 1}^J \sum_{m = 0}^M b_{jm}^\times A_{1 + m / \alpha} t^{-j - m} + O(t^{-N}).
            \label{eq: expansion of sum 1}
        \end{aligned}
        \end{equation}
        By the decay condition \eqref{eq: decay condition 1}, all the coefficients of $t^\gamma$ and $t^\gamma \log t$, $\gamma > -N$ in the expansion above (after collecting the terms) should be zero.
        
        If there exists $j \geq 1$ such that $\mu_j \neq 0$, we choose the smallest such $j$ and denote it by $j_0$. We choose $N_0$ sufficiently large so that $J \geq j_0$ when $N \geq N_0$. In this situation we extract terms of the form $t^{\gamma} \log t$ from the expansion \eqref{eq: expansion of sum 1}. Consequently, after collecting the terms, all the coefficients of $t^{-l_k \alpha - j} \log t$ in the expansion
        \begin{equation}
            \label{eq: expansion j_0 1}
            \sum_{j = j_0}^J \sum_{k = 1}^K \frac{(-1)^{l_k - 1 + j}} {\Gamma(-l_k \alpha)} \binom{-l_k \alpha - 1}{j - 1} A_{l_k} \mu_j t^{-l_k \alpha - j} \log t,
        \end{equation}
        are zero as long as $l_k \alpha + j < N$. In this sum, the exponent $-l_k \alpha - j$ of $t$ takes the maximum if and only if $j = j_0$ and $k = 1$. Since $l_1 \alpha \notin \mathbb{N} \cup \{ 0 \}$, we have $\frac{(-1)^{l_1 - 1 + j_0}} {\Gamma(-l_1 \alpha)} \neq 0$, which implies $A_{l_1} = 0$. Now the summation over $k$ in expansion \eqref{eq: expansion j_0 1} starts from $k = 2$, and the exponent $-l_k \alpha - j$ takes the maximum if and only if $j = j_0$ and $k = 2$. Repeat the process above, and we can prove that $A_{l_k} = 0$ for all $k$ such that $1 \leq k \leq K$ and $l_k \alpha + j_0 < N$. Letting $N \rightarrow \infty$, we obtain $A_{l_k} = 0$ for all $k \in \mathbb{N}$. Finally the proposition is proved by applying Lemma~\ref{lem: a_n=0} to $\{A_{l_k}\}_{k \in \mathbb{N}}$.
        
        Now suppose $\mu_j = 0$ for all $j \geq 1$, then all the coefficients $b^\times$ are zero by our convention. Hence, the expansion \eqref{eq: expansion of sum 1} reduces to
        \begin{align*}
            \sum_{n = 1}^\infty a_n \psi_n(t) 
            &= \mu_0 A_0 + \sum_{k = 1}^K \frac{(-1)^{l_k} \mu_0}{\Gamma(1 - l_k \alpha)} A_{l_k} t^{-l_k \alpha} \\
            &\quad \ + \sum_{j = 1}^J \sum_{k = 1}^K \frac{(-1)^{l_k - 1 + j}} {\Gamma(-l_k \alpha)} \binom{-l_k \alpha - 1}{j - 1} A_{l_k} t^{-l_k \alpha - j} c_{\mu, j - 1} + O(t^{-N}).
        \end{align*}
        If $\mu_0 \neq 0$, obviously $A_0$ should be zero. Then the term with largest exponent of $t$ is $\frac{(-1)^{l_1} \mu_0}{\Gamma(1 - l_1 \alpha)} A_{l_1}$ $t^{-l_1 \alpha}$. It follows that $A_{l_1} = 0$. Still we can repeat the process above to prove $A_{l_k} = 0$ for all $k \in \mathbb{N}$. When $\mu_0 = 0$, then $c_{\mu, j - 1}$ are not identically zero due to Lemma \ref{lem: mu equiv 0} and the assumption $\mu \not\equiv 0$. Using the same iteration process we can still prove $A_{l_k} = 0$ for all $k \in \mathbb{N}$. We then apply Lemma~\ref{lem: a_n=0} again to finish the proof.
    \end{pf}
    
    \begin{pf}[Proposition~\ref{prop: key prop 3}]
        Similarly to Proposition~\ref{prop: key prop 12}, we can derive an expansion for $\sum_{n = 1}^\infty \widetilde{a}_n \widetilde{\psi}_n(t)$:
        \begin{equation}
        \begin{aligned}
            \sum_{n = 1}^\infty \widetilde{a}_n \widetilde{\psi}_n(t)
            &= \widetilde{\mu}_0 \widetilde{A}_0 + \sum_{k = 1}^K \frac{(-1)^{l_k} \widetilde{\mu}_0}{\Gamma(1 - l_k \alpha)} \widetilde{A}_{l_k} t^{-l_k \alpha} \\
            & \quad \ + \sum_{j = 1}^J \sum_{k = 1}^K \frac{(-1)^{l_k - 1 + j}} {\Gamma(-l_k \alpha)} \binom{-l_k \alpha - 1}{j - 1} \widetilde{A}_{l_k} t^{-l_k \alpha - j} (\widetilde{\mu}_j \log t + c_{\widetilde{\mu}, j - 1} + \widetilde{b}_{jk}^\times) \\
            & \quad \ + \sum_{j = 1}^J \sum_{m = 0}^M \widetilde{b}_{jm}^\times \widetilde{A}_{m / \alpha} t^{-j - m} + O(t^{-N}).
            \label{eq: expansion of sum 2}
        \end{aligned}
        \end{equation}
        where
        \begin{equation*}
            \widetilde{A}_\gamma = \sum_{n = 1}^\infty \widetilde{a}_n \lambda_n^{-\gamma - 1}, \ \ \gamma \in \real.
        \end{equation*}
        By the decay condition \eqref{eq: decay condition 2}, all corresponding coefficients in the expansions \eqref{eq: expansion of sum 1} and \eqref{eq: expansion of sum 2} (after collecting the terms) should be equal.
        
        Now that $\alpha$ is irrational, so $-l_k \alpha - j = -l_{k'} \alpha - j'$ if and only if $l_k = l_{k'}$ and $j = j'$. This shows that the terms in first two sums in equation \eqref{eq: expansion of sum 1} (or \eqref{eq: expansion of sum 2}) do not overlap. A comparison of coefficients yields $A_{l_k} \mu_j = \widetilde{A}_{l_k} \widetilde{\mu}_j$ for all $j \geq 0$ and $k \in \mathbb{N}$. It follows that $\sum_{n = 1}^\infty (\mu_j a_n - \widetilde{\mu}_j \widetilde{a}_n) \lambda_n^{-l_k - 1} = 0$. By Lemma~\ref{lem: a_n=0}, we have 
        \begin{equation}
            \mu_j a_n - \widetilde{\mu}_j \widetilde{a}_n = 0
            \label{eq: mu and a}
        \end{equation}
        for all $j \geq 0$ and $n \in \mathbb{N}$.
        
        Choose $n_0 \in \mathbb{N}$ such that $\widetilde{a}_{n_0} \neq 0$. Let $a_{n_0} = \kappa \widetilde{a}_{n_0}$, then by equation \eqref{eq: mu and a} we have $\kappa \mu_j = \widetilde{\mu}_j$ for $j \geq 0$. Suppose there exists $j_0 \geq 0$ such that $\mu_{j_0} \neq 0$, then by equation \eqref{eq: mu and a} again, we get $a_n = \kappa \widetilde{a}_n$ for $n \in \mathbb{N}$. If $\widetilde{\mu} \neq \kappa \mu$, we can apply Proposition~\ref{prop: key prop 12} to $\widetilde{\mu} - \kappa \mu$ and $\{ \widetilde{a}_n \}_{n \in \mathbb{N}}$ to obtain $\widetilde{a}_n = 0$ for all $n \in \mathbb{N}$, which is a contradiction. Thus, $\widetilde{\mu} = \kappa \mu$.
        
        If $\mu_j = 0$ for all $j \geq 0$, equation \eqref{eq: mu and a} implies that $\widetilde{\mu}_j = 0$ for $j \geq 0$. Now, we can examine the expansions \eqref{eq: expansion of sum 1} and \eqref{eq: expansion of sum 2} again to deduce that $A_{l_k} c_{\mu, j} = \widetilde{A}_{l_k} c_{\widetilde{\mu}, j}$ for $j \geq 0$ and $k \in \mathbb{N}$. Another application of Lemma~\ref{lem: a_n=0} yields $c_{\mu, j} a_n - c_{\widetilde{\mu}, j} \widetilde{a}_n = 0$ for $j \geq 0$ and $n \in \mathbb{N}$. According to the assumption $\mu \not\equiv 0$ and Lemma \ref{lem: mu equiv 0}, there exists $j_0 \geq 0$ such that $c_{\mu, j_0} \neq 0$. Now we can imitate the proof of the case $\mu_{j_0} \neq 0$ to obtain $a_n = \kappa \widetilde{a}_n$ for all $n \in \mathbb{N}$ and $\widetilde{\mu} = \kappa \mu$.
    \end{pf}

    \subsection{Proof of the theorems}
    \begin{remark}
        If $f$ satisfies the regularity condition as described in Yamamoto's work [Theorem~1, \cite{y23}], our methodology for proving Theorem~\ref{thm: IP1}, based on Proposition~\ref{prop: key prop 12}, aligns closely with the approach in \cite{y23}. However, in scenarios where $f$ is characterized by a weaker regularity condition, specifically $f \in L^2(\Omega)$, a modification is necessary: in equation (3.3) of Yamamoto's demonstration, $\tilde{v} = Av$ should be substituted with $v^{(k)} = A^k v$. Despite these similarities, we opt to detail our proof to illustrate its connection to the subsequent theorems.
    \end{remark}
    
    \begin{pf}[Theorem~\ref{thm: IP1}]
        Assuming $\mu \not\equiv 0$, the function $u(x, t)$ can be represented as $u(x, t) = \sum_{n = 1}^\infty \psi_n(t) P_nf(x)$ from \eqref{eq: solution u}. To prove part (a), we arbitrarily choose a function $v \in C_c^\infty(\omega)$ and set $a_n = \langle P_nf, v \rangle_{L^2(\Omega)}$, then $\langle u(\cdot, t), v \rangle_{L^2(\Omega)} = \sum_{n = 1}^\infty a_n \psi_n(t)$. Since $||u(\cdot, t)||_{L^2(\omega)} = O(t^{-\infty}) \text{ as } t \rightarrow \infty$, we have 
        \begin{equation*}
            \left| \sum_{n = 1}^\infty a_n \psi_n(t) \right| \leq ||u(\cdot, t)||_{L^2(\omega)} ||v||_{L^2(\omega)} = O(t^{-\infty}) \text{ as } t \rightarrow \infty.
        \end{equation*}
        To apply Proposition~\ref{prop: key prop 12}, we only need to check $\sum_{n = 1}^\infty |a_n \lambda_n^{-1} \log(1 + \lambda_n)| < \infty$.
        
        For any $k \in \mathbb{N}$, let $v^{(k)} = A^k v$, then we get
        \begin{equation*}
            a_n = \lambda_n^{-k} \langle A^k P_nf, v \rangle_{L^2(\Omega)}
            = \lambda_n^{-k} \langle P_nf, A^k v \rangle_{L^2(\Omega)}
            = \lambda_n^{-k} \langle P_nf, v^{(k)} \rangle_{L^2(\Omega)}.
        \end{equation*}
        It follows that 
        \begin{equation}
            \sum_{n = 1}^\infty |a_n \lambda_n^{-1} \log(1 + \lambda_n)| 
            \leq ||f||_{L^2(\Omega)} ||v^{(k)}||_{L^2(\Omega)} \sum_{n = 1}^\infty \lambda_n^{-k - 1} \log(1 + \lambda_n).
            \label{eq: theorem 1 (a)}
        \end{equation}
        Let $\{ \widetilde{\lambda}_n \}_{n \in \mathbb{N}}$ be the eigenvalues of $A$ in ascending order and counting multiplicities. By Weyl's law, there exists a constant $c_\lambda > 0$ such that $\widetilde{\lambda}_n = c_\lambda (n^{2/d} + o(1))$. It follows that $\lambda_n \geq \frac{c_\lambda}{2} n^{2/d}$ when $n$ is large enough. By choosing a sufficiently large $k$, the series in equation \eqref{eq: theorem 1 (a)} converges. Now we can apply Proposition~\ref{prop: key prop 12} to obtain $a_n = 0$ for all $n \in \mathbb{N}$. From the choice of $v$, we get $P_nf = 0$ in $\omega$ for all $n \in \mathbb{N}$. Finally we can use the argument in Remark 2 to show that $f \equiv 0$.
        
        For part (b), we choose $v \in C_c^\infty(\gamma)$ and set $a_n = \langle \partial_\nu P_nf, v \rangle_{L^2(\partial \Omega)}$, then $\langle \partial_\nu u(\cdot, t), v \rangle_{L^2(\partial \Omega)} = \sum_{n = 1}^\infty a_n \psi_n(t)$. The condition $||\partial_\nu u(\cdot, t)||_{L^2(\gamma)} = O(t^{-\infty}) \text{ as } t \rightarrow \infty$ implies 
        \begin{equation*}
            \left| \sum_{n = 1}^\infty a_n \psi_n(t) \right| 
            \leq || \partial_\nu u(\cdot, t)||_{L^2(\gamma)} ||v||_{L^2(\gamma)} = O(t^{-\infty}) \text{ as } t \rightarrow \infty.
        \end{equation*}
        On the other hand, for any $\epsilon > 0$ we have
        \begin{equation*}
            |a_n| \leq ||\partial_\nu P_nf||_{L^2(\partial \Omega)} ||v||_{L^2(\partial \Omega)} 
            \leq C ||P_n f||_{H^{3/2 + \epsilon}(\Omega)}
            \leq C ||A^{3/4 + \epsilon / 2} P_n f||_{L^2(\Omega)}.
        \end{equation*}
        Since $f \in \mathcal{D}(A^{\sigma_2})$, we set $g_2 = A^{\sigma_2} f \in L^2(\Omega)$. Then $P_n g_2 = \lambda_n^{\sigma_2} P_n f$, and consequently
        \begin{equation*}
            |a_n| \leq C \lambda_n^{3/4 + \epsilon / 2 - \sigma_2} ||P_n g_2||_{L^2(\Omega)}.
        \end{equation*}
        Using Cauchy-Schwartz inequality, we obtain
        \begin{equation}
            \label{eq: theorem 1 (b)}
            \begin{aligned}
                \sum_{n = 1}^\infty |a_n \lambda_n^{-1} \log(1 + \lambda_n)| 
                & \leq C \left( \sum_{n = 1}^\infty ||P_n g_2||_{L^2(\Omega)}^2 \right)^{1/2} \left( \sum_{n = 1}^\infty \lambda_n^{-1/2 + \epsilon - 2\sigma_2} (\log(1 + \lambda_n))^2 \right)^{1/2} \\
                & = C ||g_2||_{L^2(\Omega)} \left( \sum_{n = 1}^\infty \lambda_n^{-1/2 + \epsilon - 2\sigma_2} (\log(1 + \lambda_n))^2 \right)^{1/2}.
            \end{aligned}
        \end{equation}
        When $n$ is large enough, we have $\lambda_n^{-1/2 + \epsilon - 2\sigma_2} (\log(1 + \lambda_n))^2 \leq \frac{c_\lambda}{2} n^{-2/d (2\sigma_2 + 1/2 - 2 \epsilon)}$. By definition of $\sigma_2$, we can choose $\epsilon$ such that $0 < \epsilon < \sigma_2 + 1/4 - d/4$. Now the series in equation \eqref{eq: theorem 1 (b)} converges, and we can argue similarly to part (a) to finish the proof.
    \end{pf}
    
    \begin{pf}[Theorem~\ref{thm: IP2}]
        We prove by contradiction that if $\mu \not\equiv 0$, then $x_0 \in \Lambda(f)$ in part (a) and $x_0 \in \Lambda_b(f)$ in part (b).
        
        For part (a) we set $a_n = P_n f(x_0)$ and for part (b) we set $a_n = \partial_\nu P_n f(x_0)$, then it is straightforward to check that $\sum_{n = 1}^\infty a_n \psi_n(t) = O(t^{-\infty})$ as $t \rightarrow \infty$. When $a_n = P_n f(x_0)$, by Sobolev embedding, for any $\epsilon > 0$ we have
        \begin{equation*}
            |a_n| \leq C ||P_nf||_{H^{d/2 + \epsilon}(\Omega)}
            \leq C ||A^{d/4 + \epsilon / 2} P_n f||_{L^2(\Omega)}.
        \end{equation*}
        Setting $g_3 = A^{\sigma_3} f \in L^2(\Omega)$, we get $|a_n| \leq C \lambda_n^{d/4 + \epsilon / 2 - \sigma_3} ||P_n g_3||_{L^2(\Omega)}$. It follows from Cauchy-Schwartz inequality that
        \begin{equation}
            \label{eq: theorem 2 (a)}
            \begin{aligned}
                \sum_{n = 1}^\infty |a_n \lambda_n^{-1} \log(1 + \lambda_n)| 
                & \leq C \left( \sum_{n = 1}^\infty ||P_n g_3||_{L^2(\Omega)}^2 \right)^{1/2} \left( \sum_{n = 1}^\infty \lambda_n^{-2 + d/2 + \epsilon - 2\sigma_3} (\log(1 + \lambda_n))^2 \right)^{1/2} \\
                & = C ||g_3||_{L^2(\Omega)} \left( \sum_{n = 1}^\infty \lambda_n^{-2 + d/2 + \epsilon - 2\sigma_3} (\log(1 + \lambda_n))^2 \right)^{1/2}.
            \end{aligned}
        \end{equation}
        When $a_n = \partial_n P_n f(x_0)$, we have 
        \begin{equation*}
            |a_n| \leq C ||\partial_\nu P_n f||_{H^{d/2 + \epsilon}(\partial \Omega)}
            \leq C ||P_nf||_{H^{(d + 3) / 2 + \epsilon}(\Omega)}
            \leq C ||A^{(d + 3) / 4 + \epsilon / 2} P_n f||_{L^2(\Omega)}.
        \end{equation*}
        Setting $g_4 = A^{\sigma_4} f \in L^2(\Omega)$, then similar to equation \eqref{eq: theorem 2 (a)} we obtain
        \begin{equation}
            \sum_{n = 1}^\infty |a_n \lambda_n^{-1} \log(1 + \lambda_n)| 
            \leq C ||g_4||_{L^2(\Omega)} \left( \sum_{n = 1}^\infty \lambda_n^{(d - 1) / 2 + \epsilon - 2\sigma_4} (\log(1 + \lambda_n))^2 \right)^{1/2}.
            \label{eq: theorem 2 (b)}
        \end{equation}
        The regularity condition of $f$ in parts (a) and (b) admits a choice of $\epsilon > 0$ such that the series in equations \eqref{eq: theorem 2 (a)} and \eqref{eq: theorem 2 (b)} converge. By Proposition~\ref{prop: key prop 12}, we have $a_n = 0$ for all $n \in \mathbb{N}$. This directly shows that $x_0 \in \Lambda(f)$ in part (a) and $x_0 \in \Lambda_b(f)$ in part (b).
    \end{pf}
    
    \begin{pf}[Theorem~\ref{thm: IP3}]
        If either $F$ or $\widetilde{F}$ is identically zero, then Proposition~\ref{prop: key prop 12} implies that the other must also be zero.
        Therefore, we focus on the scenario where neither $F$ nor $\widetilde{F}$ is zero, implying that none of $f$, $\widetilde{f}$, $\mu$ or $\widetilde{\mu}$ is identically zero.
        
        Same as in the proof of Theorem~\ref{thm: IP1}, for part (a) we choose $v \in C_c^\infty(\omega)$ and set $a_n = \langle P_nf, v \rangle_{L^2(\Omega)}$ and $\widetilde{a}_n = \langle P_n \widetilde{f}, v \rangle_{L^2(\Omega)}$. Then we have 
        \begin{equation*}
            \langle u(\cdot, t) - \widetilde{u}(\cdot, t), v \rangle_{L^2(\Omega)} 
            = \sum_{n = 1}^\infty a_n \psi_n(t) - \sum_{n = 1}^\infty \widetilde{a}_n \widetilde{\psi}_n(t).
        \end{equation*}
        For part (b) we choose $v \in C_c^\infty(\gamma)$ and set $a_n = \langle \partial_\nu P_n f, v \rangle_{L^2(\partial \Omega)}$ and $\widetilde{a}_n = \langle \partial_\nu P_n \widetilde{f}, v \rangle_{L^2(\partial \Omega)}$. Then we have 
        \begin{equation*}
            \langle \partial_\nu u(\cdot, t) - \partial_\nu \widetilde{u}(\cdot, t), v \rangle_{L^2(\partial \Omega)} 
            = \sum_{n = 1}^\infty a_n \psi_n(t) - \sum_{n = 1}^\infty \widetilde{a}_n \widetilde{\psi}_n(t).
        \end{equation*}
        Now we can proceed as in Theorem~\ref{thm: IP1} to check that the conditions in Proposition~\ref{prop: key prop 3} are satisfied, and we get $a_n = \kappa \widetilde{a}_n$ for $n \in \mathbb{N}$ and $\widetilde{\mu} = \kappa \mu$. Repeat the last part of the proof of Theorem~\ref{thm: IP1} with $f$ replaced by $f - \kappa \widetilde{f}$, and we have $f = \kappa \widetilde{f}$. Consequently, $F = f \mu = \kappa \widetilde{f} \mu = \widetilde{f} \widetilde{\mu} = \widetilde{F}$.
    \end{pf}
    
    \section{Conclusions}
    \label{sec: conclusions}
    In this paper, we have established various uniqueness results for the inverse source problem associated with the fractional diffusion-wave equation \eqref{eq: TFDE}, mainly via techniques in asymptotic analysis. Crucially, these uniqueness results depend solely on the existence of an asymptotic expansion of the temporal component $\mu$ at infinity, as opposed to any vanishing conditions on $\mu$ or the spatial component $f$.
    
    We conclude with several remarks and potential areas for future research:
    \begin{itemize}
        \item In cases where $\mu$ exhibits super-polynomial decay, condition \eqref{eq: fast decay} provides a continuous and positive function $\mu^{(1)}$ that majorizes $\mu$, and the set $\{ p \mu^{(1)} \, | \ p$ is a polynomial$\}$ is dense on $C_0[0, \infty)$. This allow us to obtain $\mu \equiv 0$ when $\int_0^\infty s^m \mu(s) \dif{s} = 0$ for all $m \geq 0$ in the proof of Proposition~\ref{prop: key prop 12}. Another possible approach is to identify a condition that ensures the denseness of $\{ p \mu \ | \ p$ is a polynomial$\}$ on $L^1(0, \infty)$ under non-vanishing $\mu$. However, a necessary condition to impose is that $\supp(\mu) = [0, \infty)$, since $p \mu$ is zero wherever $\mu$ is zero. In order to obtain the same result without the condition on the support of $\mu$, we still have to decompose $\mu$ into $\mu^{(1)} \mu^{(2)}$, where $\mu^{(1)}$ is positive and has a similar decay property to $\mu$, and $\mu^{(2)}$ takes value in $\{ -1, 0, 1 \}$. It suffices to find suitable conditions on $\mu$ or $\mu^{(1)}$ such that $\{ p \mu^{(1)} \ | \ p$ is a polynomial$\}$ is dense on $L^1(0, \infty)$. This leads to the challenge of the Bernstein approximation problem. Additionally, further conditions analogous to \eqref{eq: t_*} and \eqref{eq: fast decay} may be developed.
        
        \item The inverse problem without any vanishing assumptions and with data collected for a finite duration post-incident presents significant challenges. The current methodologies appear inadequate for this scenario. If $\mu$ is zero beyond a certain time $t_0$, the uniqueness relies on the analyticity of $\psi_n$ (defined in equation \eqref{eq: psi_n}), as discussed in \cite{y23, jk22}. Abandoning this vanishing assumption implies that $\psi_n$ only exhibits analyticity when $\mu$ itself is analytic, an impractical requirement in many cases.
            
        \item The dependency of uniqueness on the fractional order $\alpha$ is noteworthy. When $\alpha$ is irrational, the entire source term $f(x) \mu(t)$ can be determined under weaker assumptions. In contrast, when $\alpha \in (0, 1) \cup (1, 2)$, it's possible to determine only one of the two components under similar conditions. Furthermore, for $\alpha = 1$ or 2, we must impose additional conditions on $\mu$ or $f$ to achieve uniqueness, as exemplified in \cite{cly23}. For example, when $\alpha = 1$, the requirement includes $\mu \equiv 0$ on $(t_0 , \infty)$ and $\int_0^{t_0} \exp(\lambda_n s) \mu(s) \dif{s} \neq 0$ for all $n \in \mathbb{N}$. This inconsistency suggests an opportunity to explore more robust models for representing anomalous diffusion phenomena.
    \end{itemize}

    \section*{Acknowledgements}

    This work is supported by the National Key R\&D Program of China (No. 2021YFA0719200) and the National Natural Science Foundation of China (No. 11971258).

    \bibliographystyle{plainurl}
    \bibliography{ref}

\begin{thebibliography}{10}

\bibitem{cly23}
Jin Cheng, Shual Lu, and Masahiro Yamamoto.
\newblock Determination of source terms in diffusion and wave equations by
  observations after incidents: Uniqueness and stability.
\newblock {\em CSIAM Transactions on Applied Mathematics}, 4(2):381--418, 2023.
\newblock \href {https://doi.org/10.4208/csiam-am.SO-2022-0028}
  {\path{doi:10.4208/csiam-am.SO-2022-0028}}.

\bibitem{cl23}
Xing Cheng and Zhiyuan Li.
\newblock Uniqueness and stability for inverse source problem for fractional
  diffusion-wave equations.
\newblock {\em Journal of Inverse and Ill-posed Problems}, 31(6):885--904,
  2023.
\newblock \href {https://doi.org/10.1515/jiip-2021-0078}
  {\path{doi:10.1515/jiip-2021-0078}}.

\bibitem{gvh06}
E~Gerolymatou, I~Vardoulakis, and R~Hilfer.
\newblock Modelling infiltration by means of a nonlinear fractional diffusion
  model.
\newblock {\em Journal of physics D: Applied physics}, 39(18):4104, 2006.
\newblock \href {https://doi.org/10.1088/0022-3727/39/18/022}
  {\path{doi:10.1088/0022-3727/39/18/022}}.

\bibitem{gly15}
Rudolf Gorenflo, Yuri Luchko, and Masahiro Yamamoto.
\newblock Time-fractional diffusion equation in the fractional sobolev spaces.
\newblock {\em Fractional Calculus and Applied Analysis}, 18:799--820, 2015.
\newblock \href {https://doi.org/10.1515/fca-2015-0048}
  {\path{doi:10.1515/fca-2015-0048}}.

\bibitem{jk22}
Jaan Janno and Yavar Kian.
\newblock Inverse source problem with a posteriori boundary measurement for
  fractional diffusion equations.
\newblock {\em Mathematical Methods in the Applied Sciences},
  46(14):15868--15882, 2023.
\newblock \href {https://doi.org/10.1002/mma.9432}
  {\path{doi:10.1002/mma.9432}}.

\bibitem{kly22}
Yavar Kian, Yikan Liu, and Masahiro Yamamoto.
\newblock Uniqueness of inverse source problems for general evolution
  equations.
\newblock {\em Communications in Contemporary Mathematics}, 25(06):2250009,
  2023.
\newblock \href {https://doi.org/10.1142/S0219199722500092}
  {\path{doi:10.1142/S0219199722500092}}.

\bibitem{kst23}
Yavar Kian, {\'E}ric Soccorsi, and Faouzi Triki.
\newblock Logarithmic stable recovery of the source and the initial state of
  time fractional diffusion equations.
\newblock {\em SIAM Journal on Mathematical Analysis}, 55(4):3888--3902, 2023.
\newblock \href {https://doi.org/10.1137/22M1504743}
  {\path{doi:10.1137/22M1504743}}.

\bibitem{kj18}
Nataliia Kinash and Jaan Janno.
\newblock Inverse problems for a perturbed time fractional diffusion equation
  with final overdetermination.
\newblock {\em Mathematical Methods in the Applied Sciences}, 41(5):1925--1943,
  2018.
\newblock \href {https://doi.org/10.1002/mma.4719}
  {\path{doi:10.1002/mma.4719}}.

\bibitem{lly19}
Yikan Liu, Zhiyuan Li, and Masahiro Yamamoto.
\newblock Inverse problems of determining sources of the fractional partial
  differential equations.
\newblock {\em Handbook of fractional calculus with applications}, 2:411--430,
  2019.
\newblock \href {https://doi.org/10.1515/9783110571660-018}
  {\path{doi:10.1515/9783110571660-018}}.

\bibitem{lry16}
Yikan Liu, William Rundell, and Masahiro Yamamoto.
\newblock Strong maximum principle for fractional diffusion equations and an
  application to an inverse source problem.
\newblock {\em Fractional Calculus and Applied Analysis}, 19(4):888--906, 2016.
\newblock \href {https://doi.org/10.1515/fca-2016-0048}
  {\path{doi:10.1515/fca-2016-0048}}.

\bibitem{lz17}
Yikan Liu and Zhidong Zhang.
\newblock Reconstruction of the temporal component in the source term of a
  (time-fractional) diffusion equation.
\newblock {\em Journal of Physics A: Mathematical and Theoretical},
  50(30):305203, 2017.
\newblock \href {https://doi.org/10.1088/1751-8121/aa763a}
  {\path{doi:10.1088/1751-8121/aa763a}}.

\bibitem{l09}
Yury Luchko.
\newblock Maximum principle for the generalized time-fractional diffusion
  equation.
\newblock {\em Journal of Mathematical Analysis and Applications},
  351(1):218--223, 2009.
\newblock \href {https://doi.org/10.1016/j.jmaa.2008.10.018}
  {\path{doi:10.1016/j.jmaa.2008.10.018}}.

\bibitem{mk00}
Ralf Metzler and Joseph Klafter.
\newblock The random walk's guide to anomalous diffusion: a fractional dynamics
  approach.
\newblock {\em Physics reports}, 339(1):1--77, 2000.
\newblock \href {https://doi.org/10.1016/S0370-1573(00)00070-3}
  {\path{doi:10.1016/S0370-1573(00)00070-3}}.

\bibitem{n84}
RR~Nigmatullin.
\newblock To the theoretical explanation of the “universal response”.
\newblock {\em physica status solidi (b)}, 123(2):739--745, 1984.
\newblock \href {https://doi.org/10.1515/9783112495506-040}
  {\path{doi:10.1515/9783112495506-040}}.

\bibitem{n86}
RR~Nigmatullin.
\newblock The realization of the generalized transfer equation in a medium with
  fractal geometry.
\newblock {\em Physica status solidi (b)}, 133(1):425--430, 1986.
\newblock \href {https://doi.org/10.1515/9783112495483-049}
  {\path{doi:10.1515/9783112495483-049}}.

\bibitem{handbook}
Frank~W Olver, Daniel~W Lozier, Ronald~F Boisvert, and Charles~W Clark.
\newblock {\em NIST handbook of mathematical functions}.
\newblock Cambridge university press, 2010.

\bibitem{p99}
Igor Podlubny.
\newblock {\em Fractional differential equations: an introduction to fractional
  derivatives, fractional differential equations, to methods of their solution
  and some of their applications}.
\newblock elsevier, 1998.

\bibitem{p55}
Harry Pollard.
\newblock The bernstein approximation problem.
\newblock {\em Proceedings of the American Mathematical Society},
  6(3):402--411, 1955.
\newblock \href {https://doi.org/10.2307/2032781} {\path{doi:10.2307/2032781}}.

\bibitem{sy11a}
Kenichi Sakamoto and Masahiro Yamamoto.
\newblock Initial value/boundary value problems for fractional diffusion-wave
  equations and applications to some inverse problems.
\newblock {\em Journal of Mathematical Analysis and Applications},
  382(1):426--447, 2011.
\newblock \href {https://doi.org/10.1016/j.jmaa.2011.04.058}
  {\path{doi:10.1016/j.jmaa.2011.04.058}}.

\bibitem{sy11b}
Kenichi Sakamoto and Masahiro Yamamoto.
\newblock Inverse source problem with a final overdetermination for a
  fractional diffusion equation.
\newblock {\em Math. Control Relat. Fields}, 1(4):509--518, 2011.
\newblock \href {https://doi.org/10.3934/mcrf.2011.1.509}
  {\path{doi:10.3934/mcrf.2011.1.509}}.

\bibitem{s20}
Marian Slodi{\v{c}}ka.
\newblock Uniqueness for an inverse source problem of determining a
  space-dependent source in a non-autonomous time-fractional diffusion
  equation.
\newblock {\em Fractional Calculus and Applied Analysis}, 23(6):1702--1711,
  2020.
\newblock \href {https://doi.org/10.1515/fca-2020-0084}
  {\path{doi:10.1515/fca-2020-0084}}.

\bibitem{y23}
Masahiro Yamamoto.
\newblock Uniqueness for inverse source problems for fractional diffusion-wave
  equations by data during not acting time.
\newblock {\em Inverse Problems}, 39(2):024004, 2023.
\newblock \href {https://doi.org/10.1088/1361-6420/aca55c}
  {\path{doi:10.1088/1361-6420/aca55c}}.

\end{thebibliography}

\end{document}